\documentclass{elsarticle}
\def \cc#1#2#3#4{{\hbox{\sc cc}}(#1,#2,#3,#4)}
\def \tc#1#2#3#4{{\hbox{\sc tc}}(#1,#2,#3,#4)}
\def \BR#1#2{{\hbox{\sc BR}}(#1,#2)}
\def \CM#1#2{{\hbox{\sc CM}}(#1,#2)}
\def\Z{{\mathbb Z}}
\def\nm{{\,\!-\,\!}}
\def\np{{\,\!+\,\!}}
\def\ne{{\,\!=\,\!}}
\newcommand{\Sym}{\mathrm{Sym}}
\newcommand{\GG}{\mathbb{G}}
\newcommand{\HH}{\mathbb{H}}
\usepackage{amsmath}
\usepackage{amsthm}
\usepackage{amssymb}
\usepackage{xspace}
\usepackage{url}
\newcommand{\emptyword}{\epsilon}
\theoremstyle{plain}
\newtheorem{theorem}{Theorem}[section]

\newtheorem{proposition}[theorem]{Proposition}
\newtheorem{lemma}[theorem]{Lemma}

\newtheorem{remark}[theorem]{Remark}

\theoremstyle{definition}
\newtheorem{definition}[theorem]{Definition}
\newenvironment{proofof}[1]{\normalsize {\it Proof of #1}.}{{\hfill $\Box$}}
\newenvironment{mylist}{\begin{list}{}{
\setlength{\parskip}{0mm}
\setlength{\topsep}{2mm}
\setlength{\parsep}{0mm}
\setlength{\itemsep}{0.5mm}
\setlength{\labelwidth}{7mm}
\setlength{\labelsep}{3mm}
\setlength{\itemindent}{0mm}
\setlength{\leftmargin}{12mm}
\setlength{\listparindent}{6mm}
}}{\end{list}}

\begin{document}
\begin{frontmatter}
\title{Isomorphism and non-isomorphism for interval groups of type $D_n$}
	\author[bie]{Barbara Baumeister}
	\ead{baumeist@math.uni-bielefeld.de}
	\author[war]{Derek F. Holt}
	\ead{DerekHolt127@gmail.com}
	\author[bie]{Georges Neaime}
	\ead{George@Neaime@ruhr.uni-bochum.de}
	\author[ncl]{Sarah Rees\corref{cor1}}
	\ead{Sarah.Rees@newcastle.ac.uk}
	\address[bie]{Facult\"at f\"ur Mathematik, Universit\"at Bielefeld, 33615 Bielefeld, Germany}
	\address[war]{Mathematics Institute, University of Warwick,
       Coventry CV4 7AL, UK}
        \address[ncl]{School of Mathematics and Statistics, University of Newcastle,\
       Newcastle NE1 7RU, UK}
        \cortext[cor]{Corresponding author}
\begin{abstract}
	We consider presentations that were derived in \cite{BaumeisterNeaimeRees} for the interval groups associated with proper quasi-Coxeter elements of the Coxeter group $W(D_n)$. We use combinatorial methods to derive alternative presentations for the groups, and use these new presentations to show that
	the interval group associated with a proper quasi-Coxeter element of $W(D_n)$ cannot be isomorphic to the Artin group of type $D_n$.
	While the specific problems we solve
	arise from the study of interval groups, their solution provides an illustration of how techniques indicated by computational observation can be used to derive properties of
	all groups within an infinite family.
\end{abstract} 
\end{frontmatter}

\section{Introduction} 
This article uses combinatorial methods to answer a question that arises from the results of \cite{BaumeisterNeaimeRees}. 
In particular, we prove that the presentations in an infinite family associated
with the Coxeter group $W(D_n)$ define groups that are not isomorphic
to the Artin groups of the same type. Our proof of that non-isomorphism
was driven by observations made during computation with those presentations,
and provides 
an illustration of how a collection of computer results that strongly
indicate the correctness of a general conjecture for infinite classes of
groups can be used in writing a theoretical proof of those general results.

The presentations that we consider are those of interval groups
associated with quasi-Coxeter elements of the Coxeter group $W(D_n)$.
The article \cite{BaumeisterNeaimeRees}
constructs these presentations, defining each interval group
as a one-relator quotient of an Artin group $A(\Delta)$,
for a Carter diagram $\Delta$ associated with the associated
conjugacy class of quasi-Coxeter elements. 
The question of whether a group so defined is isomorphic to the Artin group
$A(D_n)$ is a natural one arising out of the theory of interval groups
associated with Coxeter (rather then quasi-Coxeter) elements of $W(D_n)$ \cite{Bessis}. 
However the question could only be partially answered using that theory
(an answer is provided for even $n$ in \cite{BaumeisterHoltNeaimeReesII}).
The proof in this article, valid for all integers $n$, uses none of that theory,
and works purely from the presentations derived in \cite{BaumeisterNeaimeRees},
using the method of Reidemeister-Schreier~\cite[Section 2.5]{HEO}.

The main results of this article 
are Theorem~\ref{thm:isom}, 
which uses combinatorial methods to find alternative presentations for the 
interval groups that we consider as one-relator
quotients of a different Artin group
(the affine Artin group $A(\tilde{A}_{n-1})$), and 
Theorem~\ref{thm:nonisom}, 
which uses those new presentations to prove that none of those interval groups
is isomorphic to the Artin group of type $D_n$ (for appropriate $n$).
Our proof of Theorem~\ref{thm:nonisom} finds subgroups of groups in the two
infinite families that would have to correspond under an isomorphism between
a pair of groups, uses Reidemeister-Schreier techniques to
compute presentations of those subgroups, and then demonstrates that these
presentations have different abelianisations, thereby establishing the
claimed non-isomorphism.

While the questions that are answered in this article
certainly arise from the study of interval groups in \cite{BaumeisterNeaimeRees},
that article provides motivation only, and no knowledge of it is necessary
for understanding of this current article. We have chosen not to define Coxeter elements or quasi-Coxeter elements of the Coxeter group $W(D_n)$ in this article, 
or to explain the construction of the interval groups associated to those 
elements, because we do not need this information; we merely study 
presentations that are found in \cite{BaumeisterNeaimeRees}.

Each of the presentations from \cite{BaumeisterNeaimeRees} that we study in this
article is associated with a {\em Carter diagram} associated with the Coxeter group $W(D_n)$.
Those Carter diagrams  are classified in \cite{Carter},
where they are called {\em admissible diagrams}. It is the diagrams
$D_n(a_i)$ ($1\leq i \leq \lfloor \frac{n}{2}\rfloor - 1$) of \cite{Carter} that are associated with the quasi-Coxeter elements of $W(D_n)$.

In this article we use the notation
$\Delta_{t,n}$ from \cite{BaumeisterNeaimeRees} for the diagram
with $n$ vertices
that is formed by attaching paths containing $t-1$ and $n-3-t$ edges to
vertices at opposite ends of a diagonal of a square,
where $1 \leq t \leq n-3$. 
For $1 \leq t \leq (n-2)/2 $, this is 
the Carter diagram $D_n(a_t$).

We use the notation $\Delta_n$ to denote the $n$-gon on $n$ vertices; when
$n$ is even this is a Carter diagram, denoted by $D_n(b_{\frac{n}{2}-1})$ in \cite{Carter}. We note that it is also the Coxeter diagram of affine type $\tilde{A}_{n-1}$

Each of the diagrams $\Delta_{t,n}$ can be found as the result of applying a
sequence of mutations to the diagram $D_n$, and hence it is covered by
\cite[Theorem 1.1]{HaleyEtAl} and \cite[Theorem A]{GM}.  
Those results prove that, for any diagram $\Gamma'$ derived from a finite type
Coxeter diagram $\Gamma$ by a sequence of mutations, the Artin group
$A(\Gamma)$ is isomorphic to the quotient of the the Artin group $A(\Gamma')$
by the normal closure of a set of {\em cycle commutators} (defined in Section~\ref{sec:notation} below)
associated with the set of chordless cycles of $\Gamma'$, each of which relates 
generators that correspond to the vertices of a chordless cycle.

Each $\Delta_{t,n}$ contains a single $4$-cycle, and $\Delta_n$ a single
$n$-cycle. So each diagram $\Delta=\Delta_{t,n}$ or $\Delta_n$ provides us with a presentation of
$A(D_n)$ as a one-relator quotient of $A(\Delta)$.
Similarly, by \cite[Theorem~3.10]{CameronEtAl} (or alternatively by later work
of \cite{BarotMarsh} that is cited by \cite{HaleyEtAl}), 
the Coxeter group $W(D_n)$ can be found as a one-relator quotient of 
$W(\Delta)$.

In Section~\ref{sec:isoms} of this paper, we re-prove some of the isomorphisms
proved as corollaries of \cite[Theorem 1.1]{HaleyEtAl}, by presenting explicit
isomorphisms, and use the same methods to prove
the following new result.  Note that we define {\em twisted cycle commutators}
in Section~\ref{sec:notation} below. See also Figure~\ref{fig:diags} below
for our numbering of the vertices of the Coxeter diagrams, and hence our 
labelling of the generators of the corresponding Artin group referred to in the theorem.

\begin{theorem}\label{thm:isom}
	For each $t$ with $1 \leq t \leq n-3 $,
the following two groups are isomorphic:
\begin{mylist}
\item[(1)] the quotient $Q=Q_{n,t}$ of $A(\Delta_{t,n})$ by the normal
closure of the twisted cycle commutator $\tc{b_1}{b_2}{b_3}{b_4}$, and
\item[(2)] the quotient $G=G_{n,t}$ of $A(\Delta_n)$ by the $t$-twisted
cycle commutator $\tc{a_1}{a_2}{\ldots}{a_n}_{t}$.
\end{mylist}
\end{theorem}
\begin{remark}
It is clear, by rotating the diagram $\Delta_{t,n}$ through $180^\circ$, that
$A(\Delta_{t,n}) \cong A(\Delta_{n-2-t,n})$, and it follows from
Lemma~\ref{lem:tcequiv} below that $Q_{n,t} \cong Q_{n,n-2-t}$ for all
$1\leq t \leq n-3$.  So we could assume that $t \leq \lfloor (n-2)/2 \rfloor $,
but this is not necessary for the proof.

It follows from the theorem that $G_{n,t} \cong G_{n,n-2-t}$, which can also
be proved directly by applying the automorphism of $A(\Delta_n)$ induced
by  $a_i \mapsto a_i^{-1}$ for $1 \leq i \leq n$.
\end{remark}

We might expect the groups $Q$ and $G$ to have soluble word problems,
but we have been unable to prove that.

In Section~\ref{sec:nonisoms} we prove Theorem~\ref{thm:nonisom} (stated below) that
the groups $Q_{n,t}$ and $G_{n,t}$ are not isomorphic to the Artin group of type $D_n$.
This theorem contrasts with the results of \cite{HaleyEtAl} that we have already mentioned,
and which we re-prove in Section 3 as
Propositions~\ref{PropArtin-n-gon},~\ref{PropArtin-Flag}~and~\ref{Prop-RattleEquiv}, which find isomorphisms between
$A(D_n)$ and each of the quotients 
$A(\Delta_{t,n})$ by the normal closure of the cycle commutator
$\cc{b_1}{b_2}{b_3}{b_4}$, and of $A(\Delta_n)$ by the normal
closure of the cycle commutator $\cc{a_1}{a_2}{\ldots}{a_n}$.
We note that each of these presentations for quotients by cycle commutators can be derived in a natural 
way as a presentation for the interval group of a Coxeter (rather than quasi-Coxeter) element of type 
$D_n$, on a subset of the set of all reflections (using methods described in \cite{BaumeisterNeaimeRees}). 

\begin{theorem}\label{thm:nonisom}
The groups $Q_{n,t}$ and $G_{n,t}$ referred to in Theorem~\ref{thm:isom}
are not isomorphic to $A(D_n)$ for any $n,t$ with $n \geq 4$
and $1 \leq t \leq n-3$.
\end{theorem}

This result for $n$ up to about $70$ was originally proved by computer
calculations, in which we found subgroups of the groups $G$ of
Theorem~\ref{thm:isom} and of $A(D_n)$ that would have to correspond under
a putative isomorphism, but which had different abelianisations. We were
able to observe the details of the steps in the computer calculations, and
to use these to construct a proof for general $n$.

We note that the article \cite{BaumeisterHoltNeaimeReesII} contains a proof that 
the groups $Q_{n,t}$ are non-isomorphic to $A(D_n)$ when $n$ is even, but
that proof does not extend to the case where $n$ is odd. 
We have not resolved the question of whether, for a given value of $n$, the
groups $G_{n,t}$ are isomorphic for different values of $t$, except
for the isomorphism $G_{n,t} \cong G_{n,n-2-t}$ mentioned in the remark above.

\textbf{Acknowledgements}. The third author would like to acknowledge support through the DFG grants BA2200/5-1 and RO1072/19-1. The first three authors would like to thank Sarah Rees who hosted them during their research stay in Newcastle.

All four authors would like to acknowledge some very helpful and constructive comments from the referees.

\section{Notation and preliminary lemmas}\label{sec:notation}

Throughout this article we denote the identity element of a group by $1$. 
Given words $u,v$ we write $u=v$ to indicate that words are identical as strings, and $u=_G v$ to indicate that $u,v$ represent the same element of the group $G$.

\begin{figure}
\mbox{
{\small
\begin{picture}(120,130)(-60,-10)
\put(-5,-10){$\Delta_n$}
\put(-10,95){\circle*{4}}
\put(-10,95){\line(1,0){20}}
\put(-14,100){$n$}
\put(10,95){\circle*{4}}
\put(8,100){$1$}

\put(-10,95){\line(-2,-1){18}}
\put(-28,86){\circle*{4}}
\put(-52,84){$n\nm 1$}

\put(10,95){\line(2,-1){18}}
\put(28,86){\circle*{4}}
\put(32,84){$2$}

\put(-28,86){\line(-1,-2){8.5}}
\put(-36.5,69){\circle*{4}}
\put(-61,66){$n\nm 2$}

\put(28,86){\line(1,-2){8.5}}
\put(36.5,69){\circle*{4}}
\put(40.5,66){$3$}

\put(-36.5,49){\circle*{1.8}}
\put(-28,29){\circle*{1.8}}
\put(-19,25.5){\circle*{1.8}}
\put(-10,22){\circle*{1.8}}
\put(-3,20){\circle*{1.8}}
\put(36.5,49){\circle*{1.8}}
\put(28,29){\circle*{1.8}}
\put(19,25.5){\circle*{1.8}}
\put(10,22){\circle*{1.8}}
\put(3,20){\circle*{1.8}}
\end{picture}
\begin{picture}(120,130)(-50,-70)
\put(-5,-70){$\Delta_{t,n}$}
\put(-8,-8){\circle*{4}}
\put(-18,-12){$4$}
\put(-8,8){\circle*{4}}
\put(-18,3){$1$}
\put(8,-8){\circle*{4}}
\put(12,-10){$3$}
\put(8,8){\circle*{4}}
\put(12,5){$2$}
\put(-8,-8){\line(1,0){16}}
\put(-8,8){\line(1,0){16}}
\put(-8,-8){\line(0,1){16}}
\put(8,-8){\line(0,1){16}}

\put(-8,8){\line(-1,1){11}}
\put(-19,19){\circle*{4}}
\put(-15,18){$5$}
\put(-19,19){\line(-1,1){11}}
\put(-30,30){\circle*{4}}
\put(-26,29){$6$}
\put(-45,45){\circle*{4}}
\put(-41,45){$n \np 1 \nm t \ne r+4$}

\put(8,-8){\line(1,-1){11}}
\put(19,-19){\circle*{4}}
\put(-3,-24){$r \np 5$}
\put(19,-19){\line(1,-1){11}}
\put(30,-30){\circle*{4}}
\put(8,-35){$r \np 6$}
\put(45,-45){\circle*{4}}
\put(-4,-49){$r \np s \np 4 \ne n$}

\put(-35,35){\circle*{1.8}}
\put(35,-35){\circle*{1.8}}
\put(-40,40){\circle*{1.8}}
\put(40,-40){\circle*{1.8}}
\end{picture}
\begin{picture}(120,130)(-10,-70)
\put(35,-70){$D_n$}
\put(-14,-14){\circle*{4}}
\put(-10,-20){$2$}
\put(-14,14){\circle*{4}}
\put(-10,16){$1$}
\put(0,0){\circle*{4}}
\put(2,4){$3$}
\put(-14,-14){\line(1,1){14}}
\put(-14,14){\line(1,-1){14}}
\put(0,0){\line(1,0){20}}
\put(20,0){\circle*{4}}
\put(20,4){$4$}
\put(20,0){\line(1,0){20}}
\put(40,0){\circle*{4}}
\put(40,4){$5$}

\put(75,0){\circle*{4}}
\put(67,4){$n \nm 1$}
\put(75,0){\line(1,0){20}}
\put(95,0){\circle*{4}}
\put(96,4){$n$}

\put(47,0){\circle*{1.8}}
\put(54,0){\circle*{1.8}}
\put(61,0){\circle*{1.8}}
\put(68,0){\circle*{1.8}}
\end{picture}
} }

\caption{The labelled diagrams}\label{fig:diags}
\end{figure}
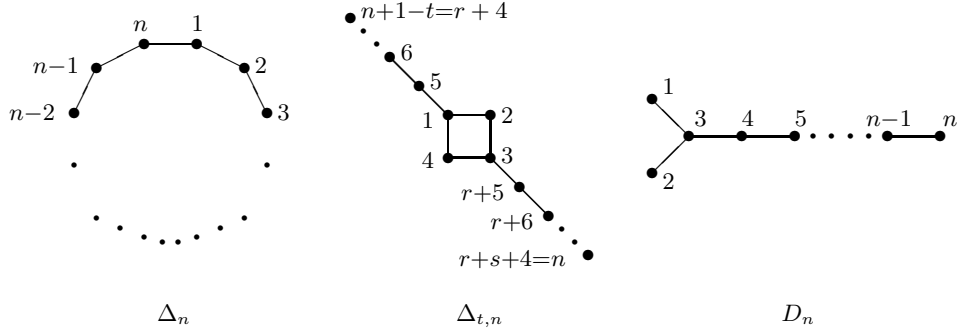

The labelling of the vertices in the diagrams that we shall now describe
can be seen in Figure~\ref{fig:diags}.

We label the vertices of the Coxeter diagram of $D_n$ as $1,2,\ldots,n$,
with $1$ and $2$ labelling the two vertices at ends of the fork, 
$3$ joined to both $1$ and $2$, and then $4,5,6,\ldots,n$ labelling the 
successive vertices from $3$ to the end of the diagram.  
We label the generators of the Artin group $A(D_n)$ as $x_1,\dots,x_n$.
So $x_1$ and $x_2$ commute, $x_1$ and $x_2$ braid with $x_3$,
the generator $x_i$ braids with $x_{i+1}$ for $3 \leq i \leq n-1$,
and all other pairs of generators commute.

We label the vertices of the $\Delta_n$ diagram
(this is $\tilde{A}_{n-1}$) as $1,2,\ldots, n$ around the cycle.
We label the generators of the Artin group $A(\Delta_n)$ as
$a_1,\dots,a_n$. 
So, for $1\leq i \leq n-1$, $a_i$ braids with $a_{i+1}$, $a_n$ braids with $a_1$ and all other pairs of 
generators commute.

We label the vertices of the $\Delta_{1,n}$ diagram
$1,2,\ldots, n$ with $1,2,3,4$ around a square, and $1,5,6,\ldots,n$ labelling
consecutive vertices on a path.
We label the generators of the Artin group $A(\Delta_{1,n})$ as
$b_1,\dots,b_n$.  So the pairs $\{b_1,b_2\}$, $\{b_2,b_3\}$, $\{b_3,b_4\}$,
$\{b_4,b_1\}$, $\{b_1,b_5\}$ and then $\{b_i,b_{i+1}\}$ for $5 \leq i \leq n-1$
braid, and all other pairs of generators commute.

Define $r:=n-3-t$, and $s:=t-1$, so that the two arms of $\Delta_{t,n}$ contain
$r$ and $s$ edges respectively. Some of our constructions may seem clearer when
expressed in terms of the parameters $r,s$ rather than $t$.
For $t>1$, we label the vertices of the $\Delta_{t,n}$ diagram
$1,2,\ldots, n$, with $1,2,3,4$ around a square,
and $1,5,6,\ldots,4+r$ and $3,5+r,\ldots,n=4+r+s$ each labelling
sequential vertices on a path (where we recall that $r=n-3-t$ and $s=t-1$).
We label the generators of the Artin group $A(\Delta_{t,n})$ as
$b_1,\dots,b_{4+r}$, $c_{5+r},\ldots,c_{n}$.
So the pairs $\{b_1,b_2\}$, $\{b_2,b_3\}$, $\{b_3,b_4\}$, $\{b_4,b_1\}$;
$\{b_1,b_5\}$ and $\{b_i,b_{i+1}\}$ for $5 \leq i \leq 3+r$;
and $\{b_3,c_{5+r}\}$ and $\{c_i,c_{i+1}\}$ for $5+r \leq i \leq n-1$
braid, while all other pairs of vertices commute.

\begin{definition}
For $n \geq 3$, we define the cycle commutator $\cc{y_1}{y_2}{\dots}{y_n}$
to be the commutator
$[y_1,y_2\ldots y_{n-1}y_ny_{n-1}^{-1}\ldots y_2^{-1}]$.
So we have, for example, 
$\cc{y_1}{y_2}{y_3}{y_4} = [y_1,y_2y_3y_4y_3^{-1}y_2^{-1}]$.

We define the twisted cycle commutator
$\tc{y_1}{y_2}{y_3}{y_4}$ to be the commutator
$[y_1,y_2^{-1}y_3y_4y_3^{-1}y_2]$ and, for $n \geq 4$ and $1 \leq t \leq n-2$,
we define the $t$-twisted cycle commutator
${\tc{y_1}{y_2}{\dots}{y_n}}_t$ to be the commutator
\[[y_1,y_2^{-1}\cdots y_{t+1}^{-1}y_{t+2}\cdots
   y_{n-1}y_ny_{n-1}^{-1}\cdots y_{t+2}^{-1}y_{t+1}\cdots y_2].\]
So in particular $\tc{y_1}{y_2}{y_3}{y_4}_1 = \tc{y_1}{y_2}{y_3}{y_4},$
and ${\tc{y_1}{y_2}{\dots}{y_n}}_{n-2}$ is the commutator
\[[y_1,y_2^{-1}\cdots y_{n-1}^{-1}y_ny_{n-1}\cdots y_2].\]
\end{definition}

In the above definitions, the reader should feel free to use whichever of the
two possible definitions of the commutator $[g,h]$ they prefer -
it will make no difference! We shall show in Lemmas~\ref{lem:ccequiv}
and~\ref{lem:tcequiv} below that, in the context of our applications, these
definitions have various equivalent formulations.

For group elements $g,h$, we will write $\CM{g}{h}$ to mean that $g$ and $h$
commute (i.e.\ $gh=_G hg$) and $\BR{g}{h}$ to mean that $g$ and $h$ braid
(i.e.\ $ghg=_G hgh$.) 
We will sometimes write statements such as $\BR{g}{h} = \BR{g'}{h'}$.
By this we mean that the two statements $\BR{g}{h}$ and $\BR{g'}{h'}$ are logically equivalent.
We will sometimes refer to the {\em conjugate} of such a relation, by an element $k$; by this we mean a relation of the form
$k^{-1}ghk =_G k^{-1}hgk$ or $k^{-1}ghgk =_G k^{-1}hghk$.

The following technical lemmas will be used repeatedly in Section~\ref{sec:isoms}
\begin{lemma}\label{lem:brcm}
Let $G$ be a group, and suppose that $f,g,h \in G$ satisfy
$\BR{f}{g}$, $\BR{g}{h}$ and $\CM{f}{h}$. Then
\begin{mylist}
\item[(i)] $\BR{f}{ghg^{-1}}$,  $\BR{f}{g^{-1}hg}$, and
  $\BR{f}{fghg^{-1}f^{-1}}$;
\item[(ii)] $\CM{g}{fghg^{-1}f^{-1}}$, $\CM{f^{-1}gf}{h^{-1}gh}$, and
$\CM{fgf^{-1}}{hgh^{-1}}$.
\end{mylist}
\end{lemma}
\begin{proof}  Since $hghg^{-1}h^{-1} =_G g$ and $h$ commutes with $f$,
we have
\[ \BR{f}{g} \Rightarrow \BR{f}{hghg^{-1}h^{-1}}
\Rightarrow \BR{f}{ghg^{-1}}
\]
and similarly $h^{-1}g^{-1}hgh =_G g$ implies $\BR{f}{g^{-1}hg}.$
Then $\BR{f}{fghg^{-1}f^{-1}}$ follows by conjugating
$\BR{f}{ghg^{-1}}$ by $f$.

Furthermore, putting $w = fghg^{-1}f^{-1}$, we have
\[ gw=gfghg^{-1}f^{-1} =_G fgfhg^{-1}f^{-1} =_G fghfg^{-1}f^{-1} =_G
fghg^{-1}f^{-1}g=wg, \] 
proving $\CM{g}{fghg^{-1}f^{-1}}$.
Also, conjugating $\CM{f}{h}$, we derive both $\CM{g^{-1}fg}{g^{-1}hg} =
\CM{fgf^{-1}}{hgh^{-1}}$ and $\CM{gfg^{-1}}{ghg^{-1}} =\linebreak
\CM{f^{-1}gf}{h^{-1}gh}$.
\end{proof}

\begin{lemma}\label{lem:longconj} Let $G$ be a group, and suppose that $y_1,\ldots,y_k \in G$ (with $k \geq 2$)
satisfy $\BR{y_i}{y_{i+1}}$ for $1 \leq i < k$ and
$\CM{y_i}{y_j}$ for $1 \leq i < j-1 \leq k-1$ (as in the braid group).
Let $w := y_1 y_2 \cdots y_{k-1} y_k y_{k-1}^{-1} \cdots y_2^{-1}y_1^{-1}$.
Then
\begin{mylist}
\item[(i)] $ w =_G
y_k^{-1} y_{k-1}^{-1} \cdots y_2^{-1} y_1 y_2 \cdots y_{k-1} y_k;$
\item[(ii)] $\BR{y_1}{w}$, $\BR{y_k}{w}$, and (if $k\geq 3$) $\BR{y_k}{v}$,\\
	where $v :=  y_{k-1}^{-1} \cdots y_2^{-1} y_1 y_2 \cdots y_{k-1}$;
\item[(iii)] $\CM{y_i}{w}$ for $2 \leq i \leq k-1$. 
\end{mylist}
\end{lemma}
\begin{proof} (i) We prove this by induction on $k$. It follows
directly from $\BR{y_1}{y_2}$ when $k=2$. For $k>2$, we use
$\BR{y_{k-1}}{y_k}$ to replace the subword $y_{k-1} y_k y_{k-1}^{-1}$
of $w$ by $y_k^{-1}y_{k-1}y_k$, and then use
$\CM{y_i}{y_k}$ for $1 \leq i \leq k-2$ and induction to complete the proof.

(ii) To get $\BR{y_k}{w}$, note that $\BR{y_k}{y_{k-1}}$ implies
	$\BR{y_k}{y_k^{-1}y_{k-1}y_k} = \BR{y_k}{y_{k-1}y_ky_{k-1}^{-1}}$, and
conjugating this by $y_1 \cdots y_{k-2}$, which commutes with $y_k$,
gives $\BR{y_k}{w}$. This implies $\BR{y_k}{y_kwy_k^{-1}}= \BR{y_k}{v}$.
Similarly, it follows from $\BR{y_1}{y_2}$ that $\BR{y_1}{y_1y_2y_1^{-1}} =
\BR{y_1}{y_2^{-1}y_1y_2}$ and (when $k \geq 3$) conjugating by
$y_k^{-1} \cdots y_3^{-1}$  gives
$\BR{y_1}{y_k^{-1}\cdots y_2^{-1}y_1
y_2 \cdots y_k}$, and then $\BR{y_1}{w}$ follows from Part (i).
	When $k=2$, $\BR{y_1}{w}$ follows from $\BR{y_1}{y_2}$.

(iii) For a given $i$ with $2 \leq i \leq k-1$, let $f=y_{i-1}$, $g=y_i$, and
$h=y_{i+1} \cdots y_{k-1}y_k y_{k-1}^{-1} \cdots y_{i+1}^{-1}$.
Then we have $\CM{f}{h}$ and $\BR{f}{g}$, and $\BR{g}{h}$ follows from (ii), so
$\CM{y_i}{y_{i-1}y_i \cdots y_{k-1}y_k y_{k-1}^{-1} \cdots
y_i^{-1}y_{i-1}^{-1}}$ follows from Lemma~\ref{lem:brcm}\,(ii).
Conjugating by $y_1 \cdots y_{i-2}$ gives $\CM{y_i}{w}$.
\end{proof}

The next lemma is also proved as \cite[Lemma 2.4]{GM}. We repeat the proof
here because we will need to use essentially the same argument in the
following lemma on twisted cycle commutators.

\begin{lemma}\label{lem:ccequiv}
Let $n \geq 3$, and suppose that $a_1,a_2,\ldots,a_n$ are elements of a
group $G$ that satisfy the defining relations of $A(\Delta_n)$.
Then the $n$ properties
$$\cc{a_1}{a_2}{\dots}{a_n}=_G 1,\, \cc{a_2}{a_3}{\dots}{a_1} =_G 1,\ldots,
\cc{a_n}{a_1}{\dots}{a_{n-1}} =_G 1$$ are all equivalent.
\end{lemma}

\begin{proof}
We have
\begin{eqnarray*}
\cc{a_1}{a_2}{\dots}{a_n}=_G 1&\Rightarrow&\\
\CM{a_1}{a_2\ldots a_{n-1}a_na_{n-1}^{-1}\ldots a_2^{-1}}&\Rightarrow&\\
\CM{a_na_1a_n^{-1}}{a_na_2\ldots a_{n-1}a_na_{n-1}^{-1}\ldots
   a_2^{-1}a_n^{-1}} &\Rightarrow&\\
\CM{a_na_1a_n^{-1}}{a_2\ldots
   a_{n-2}a_na_{n-1}a_na_{n-1}^{-1}a_n^{-1}a_{n-2}^{-1}\ldots a_2^{-1}}&=_G&\\
\CM{a_na_1a_n^{-1}}{a_2\ldots a_{n-2}a_{n-1}a_{n-2}^{-1}\ldots a_2^{-1}}
&\Rightarrow&\\
\CM{a_1a_na_1a_n^{-1}a_1^{-1}}{a_1a_2\ldots
  a_{n-2}a_{n-1}a_{n-2}^{-1}\ldots a_2^{-1}a_1^{-1}} &=_G&\\
\CM{a_n}{a_1a_2\ldots a_{n-2}a_{n-1}a_{n-2}^{-1}\ldots a_2^{-1}a_1^{-1}}
\end{eqnarray*}
and so $\cc{a_n}{a_1}{\dots}{a_{n-1}} =_G 1$, and we derive the equivalence
of the properties by iterating this argument.
\end{proof}

\begin{lemma}\label{lem:tcequiv}
Let $n \geq 3$ and $1 \leq t \leq n-2$, and suppose that $a_1,a_2,\ldots,a_n$
are elements of a
group $G$ that satisfy the defining relations of $A(\Delta_n)$. 
Then the following statements hold.
\begin{mylist}
\item[(i)]
The $n$ properties
\begin{eqnarray*}
\tc{a_1}{a_2}{\dots}{a_n}_t&=_G& 1,\,
\tc{a_2}{a_3}{\dots}{a_1}_t =_G 1,\ldots,\\
\tc{a_n}{a_1}{\dots}{a_{n-1}}_t &=_G& 1
\end{eqnarray*}
are all equivalent.
\item[(ii)] Let $I := \{2,3,\ldots,n-1\}$ and let $S$ be a subset of
$I$ with $|S|=t$. Define $\epsilon_i$ for $i \in I$ by
$\epsilon_i= -1$ if $i \in S$ and $\epsilon_i = 1$ if $i \not\in S$.
\[\tc{a_1}{a_2}{\dots}{a_n}_t =_G 1 \Longleftrightarrow
[a_1, a_2^{\epsilon_2} \cdots a_{n-1}^{\epsilon_{n-1}} a_n
  a_{n-1}^{-\epsilon_{n-1}} \cdots  a_2^{-\epsilon_2}] =_G 1.\]
\end{mylist}
\end{lemma}
\begin{proof}
We prove (ii) first. There is nothing to prove if $t = n-2$, so suppose that
$t < n-2$ (and hence $n \geq 4$). Let $k \in I$ with $k < n-2$,
and for $i \in I \setminus\{k,k+1\}$, choose $\epsilon_i = \pm 1$. Define
\begin{eqnarray*}
c_1 &:=& a_2^{\epsilon_2} \cdots a_{k-1}^{\epsilon_{k-1}}
  a_k^{-1} a_{k+1} a_{k+2}^{\epsilon_{k+2}} \cdots a_{n-1}^{\epsilon_{n-1}},\\
c_2 &:=& a_2^{\epsilon_2} \cdots a_{k-1}^{\epsilon_{k-1}}
  a_k a_{k+1}^{-1} a_{k+2}^{\epsilon_{k+2}} \cdots a_{n-1}^{\epsilon_{n-1}}.
\end{eqnarray*}
We shall show below that, for any choice of $k$ and $\epsilon_i$,
we have
\[ [a_1,c_1 a_n c_1^{-1}]=_G 1 \Longleftrightarrow
	[a_1,c_2 a_n c_2^{-1}] =_G 1.\quad (*) \]
Once $(*)$ is proved, (ii) will follow. For suppose $S$ of size $t$ is given. There must exist $i,i+1 \in I$ such that exactly one of $i$ and $i+1$ is 
	in $S$, and then the set $S'$ formed from $S$ by deleting from $S$ one of those elements of $\{i,i+1\}$ and adding in the other one, also has size $t$;
	a sequence of such `moves' can be found to define a sequence of sets
	$S=S_0,S_1,\dots S_n=\{2,3,\ldots,t+1\}$, for which any two sets that are adjacent in the sequence differ
	by a single pair of elements $i,i+1$ of $I$. 
	We associate the commutator
$[a_1, a_2^{\epsilon_2} \cdots a_{n-1}^{\epsilon_{n-1}} a_n
  a_{n-1}^{-\epsilon_{n-1}} \cdots  a_2^{-\epsilon_2}]$ with the set $S$
  and 
	$\tc{a_1}{a_2}{\dots}{a_n}_t$ with the set $\{2,\ldots,t+1\}$.
	We use $(*)$ to prove equivalence
	of the relations defined by commutators associated with successive terms in that sequence, and then (ii) follows.

	To avoid horrific notation, we shall write out the proof of $(*)$ in the specific case
$n = 6$, $k=3$, $t=2$, $\epsilon_2=1$, $\epsilon_5 = -1$,
which we hope will make
the general case clear.
In this case $c_1 = a_2 a_3^{-1} a_4 a_5^{-1}$ and $c_2 = a_2 a_3 a_4^{-1} a_5^{-1}$.
Using the braid and commutativity relations in $G$,
we have
\begin{eqnarray*}
[a_1,c_1 a_n c_1^{-1}] =
[a_1,a_2 a_3^{-1} a_4 a_5^{-1}a_6 a_5 a_4^{-1} a_3 a_2^{-1}] =_G 1
                        & \Longleftrightarrow &\\
{}[a_1,a_4^{-1}a_2 a_3^{-1} a_4 a_5^{-1}a_6 a_5 a_4^{-1}a_3 a_2^{-1}a_4] =_G 1
                       & \Longleftrightarrow &\\
{}[a_1,a_2a_4^{-1}a_3^{-1} a_4 a_5^{-1}a_6 a_5 a_4^{-1}a_3a_4a_2^{-1}] =_G 1
                       & \Longleftrightarrow &\\
{}[a_1,a_2a_3 a_4^{-1}a_3^{-1}a_5^{-1}a_6 a_5 a_3a_4a_3^{-1}a_2^{-1}] =_G 1
                       & \Longleftrightarrow &\\
{}[a_1,a_2a_3 a_4^{-1}a_5^{-1}a_6 a_5 a_4a_3^{-1}a_2^{-1}]
  = [a_1,c_2 a_n c_2^{-1}] =_G 1.
\end{eqnarray*}

For (i) with $t<n-2$, we can prove as in the proof of Lemma~\ref{lem:ccequiv}
that $\tc{a_1}{a_2}{\dots}{a_n}_t =_G 1$ if and only if
\[\CM{a_n}{a_1a_2^{-1}\cdots a_{t+1}^{-1}a_{t+2} \cdots
  a_{n-2}a_{n-1}a_{n-2}^{-1}\cdots a_{t+2}^{-1}a_{t+1} \cdots a_2a_1^{-1}},\]
and the result now follows from (ii). The proof when $t=n-2$ is similar,
but in that case we conjugate the commuting relation by $a_n^{-1}$ and then by $a_1^{-1}$
rather than by $a_n$ and $a_1$.
\end{proof}

\section{Specific transformations to show isomorphisms between $A(D_n)$
and the quotients of $A(\Delta)$}\label{sec:isoms}

Although Propositions~\ref{PropArtin-n-gon}, \ref{PropArtin-Flag}, and
\ref{Prop-RattleEquiv} can be found as corollaries of \cite[Theorem 1.1]{HaleyEtAl}, we have provided our own proofs, in order to see explicit isomorphisms.
Our construction of the proofs of those Propositions also helped us to construct
the proof of Theorem~\ref{thm:isom}, which is a new result.

Our proofs of these three propositions and of Theorem~\ref{thm:isom} below
all have the following structure. We are trying to establish that
two groups $A$ and $B$ defined by presentations on generators
$\alpha_1,\ldots,\alpha_n$ and $\beta_1,\ldots,\beta_n$ are isomorphic.
To do this, we first define homomorphisms $\phi:A \to B$ and $\psi:B \to A$,
and then prove that $\phi$ and $\psi$ are mutually inverse, which implies
that they are isomorphisms. We construct $\phi$ by first defining
$\phi(\alpha_1),\ldots,\phi(\alpha_n)$ as words over the $\beta_i$
(and we abuse notation by denoting these words by
$\alpha_1,\ldots,\alpha_n$) and then verify that these
images satisfy the defining relations of the group $A$, which proves that
the map $\phi$ defined on the generators extends to a homomorphism $A \to B$.
The homomorphism $\psi:B \to A$ is defined similarly. To prove that $\phi$ and
$\psi$ are mutually inverse, it is sufficient to verify this on the
generators; that is, we show that $\psi(\phi(\alpha_i)) = \alpha_i$ and
$\phi(\psi(\beta_i)) = \beta_i$ for $1 \le i \le n$.

\begin{proposition}\label{PropArtin-n-gon}
If $n\geq 4$, then the Artin group $G:=A(D_n)$ is isomorphic to the quotient
$Q$ of the Artin group $A(\Delta_n)$ on generators
$a_1,\ldots,a_n$ by the normal closure of the cycle commutator
$\cc{a_1}{a_2}{\dotsc}{a_n}$.
\end{proposition}

\begin{proof}
Where $x_1,\ldots,x_n$ are generators of the group  $G$ 
in its standard presentation, we define
elements $a_1,\ldots,a_n$ of $G$ as words over $x_1,\ldots,x_n$ as follows.
	\begin{eqnarray*}
&& a_1 := x_1,\,a_2 := x_3,\,a_3 := x_4,\ldots, a_{n-1} := x_n,\\
&&a_n := x_n^{-1}x_{n-1}^{-1}\cdots x_3^{-1}x_2x_3\cdots x_{n-1}x_n.\end{eqnarray*}
Given the Artin relations between the $x_i$, 
we need to verify that all of the relations of $A(\Delta_n)$ hold between
the $a_i$, and also verify the cycle commutator relation
$\cc{a_1}{a_2}{\dotsc}{a_n} =_G 1$.
That will prove that the  map from $Q$ to $G$ defined by mapping
each generator $a_i$ of $Q$ to the element $a_i$ of $G$ that we have just defined
extends to a homomorphism from $Q$ to $G$.

We have
\[a_n = a_{n-1}^{-1}a_{n-2}^{-1}\cdots a_2^{-1}x_2a_2\cdots a_{n-2}a_{n-1}
\ \mbox{so}\ x_2 =a_2a_3 \cdots a_{n-1}a_na_{n-1}^{-1} \cdots a_2^{-1},\]
and $\cc{a_1}{a_2}{\dotsc}{a_n} =_G 1$ follows immediately from $\CM{x_1}{x_2}$.

Note that, since $x_4,\ldots,x_n$ commute with $x_1$, $\BR{a_1}{a_n}$
is equivalent to $\BR{x_1}{x_3^{-1}x_2x_3}$,
	and that this (conjugating) is equivalent to $\BR{x_2}{x_3x_1x_3^{-1}},$ which follows from
Lemma~\ref{lem:brcm}\,(i) applied with $f=x_2$, $g=x_3$, $h=x_1$.

The other relations are immediate apart from those that involve
$a_n$: i.e. $\BR{a_{n-1}}{a_n}$ and $\CM{a_i}{a_n}$ for $2 \leq i \leq n-2$.  
They all follow from Lemma~\ref{lem:longconj}\,(ii) and (iii).

Conversely, given generators $a_1,\ldots,a_n$ of $Q$, we define
elements $x_1,\ldots,x_n$ of $Q$ as products of those.
\[ x_1:= a_1,\, 
x_2 := a_2a_3\cdots a_{n-1}a_na_{n-1}^{-1}\cdots a_3^{-1}a_2^{-1},
x_3 := a_2,\ldots,x_n:= a_{n-1}\]
Then we need to verify that all the relations of $G$ hold between the $x_i$.
That will prove that the map from $G$ to $Q$ defined by mapping each generator
$x_i$ of $G$ to the element $x_i$ of $Q$ that we have just defined
extends to a homomorphism from $G$ to $Q$.

The relation $\CM{x_1}{x_2}$ is precisely $\cc{a_1}{a_2}\cdots{a_n}=_Q 1$.
All other relations not involving $x_2$ follow immediately from those of
the $n$-gon. So that leaves $\BR{x_2}{x_3}$ and $\CM{x_2}{x_k}$ for
$4 \leq k \leq n$. These all follow from Lemma~\ref{lem:longconj}\,(ii) and (iii).

We also need to check that the two maps between the two groups are mutually
inverse, but this is clear from their actions on the generators.
\end{proof}

\begin{proposition}\label{PropArtin-Flag}
If $n\geq 4$, then the Artin group $G := A(D_n)$ is isomorphic to the
quotient $Q$ of the Artin group $A(\Delta_{1,n})$ on generators
$b_1,\ldots,b_n$ by the normal closure of the cycle commutator
$\cc{b_1}{b_2}{b_3}{b_4}$.
\end{proposition}

\begin{proof}
{\ }
Where $x_1,\ldots,x_n$ are the generators of $G$
in its standard presentation we define
elements $b_1,\ldots,b_n$ of $G$ as products of those as follows.
\begin{eqnarray*}
&&b_1 := x_4,\, b_2 := x_3,\, b_3 := x_1,\,
b_4 := x_1^{-1}x_3^{-1}x_2x_3x_1,\\
&&b_k := x_k\ \mbox{for}\ 5 \leq k \leq n.
\end{eqnarray*}
We need now to check that the relations of $Q$ hold between these elements of $G$.
All relations of $Q$ involving $b_k$ for $k\geq 5$ are immediate, as are
$\BR{b_1}{b_2}$, $\BR{b_2}{b_3}$ and $\CM{b_1}{b_3}$ whereas, since
$b_4$ freely reduces to $x_2x_3x_1x_3^{-1}x_2^{-1}$ by Lemma~\ref{lem:longconj}\,(i),
$\BR{b_3}{b_4}$ and $\CM{b_2}{b_4}$ come from Lemma~\ref{lem:longconj}\,(ii)
and~(iii).  Furthermore, since $b_4 =_G x_2x_1^{-1}x_3x_1x_2^{-1}$ and
$x_1$ and $x_2$ commute with $x_4$, we see that $\BR{b_4}{b_1}$ follows from
$\BR{x_3}{x_4}$, whereas, since 
$b_2b_3b_4b_3^{-1}b_2^{-1} =_G x_2$,
	$\cc{b_1}{b_2}{b_3}{b_4}=_G 1$ follows immediately from
$\CM{x_2}{x_4}$. 

In the other direction, we define
elements $x_i$ of $Q$ in terms of the generators $b_i$ of $Q$ as follows.
\begin{eqnarray*}
&&x_1 := b_3,\, x_2 := b_2b_3b_4b_3^{-1}b_2^{-1},\,
x_3 := b_2,\, x_4 := b_1,\\
&&x_k := b_k\ \mbox{for}\ 5 \leq k \leq n.
\end{eqnarray*}
All relations of $G$ involving $x_k$ for $k\geq 5$ are immediate, as are
$\BR{x_1}{x_3}$, $\BR{x_4}{x_3}$ and $\CM{x_1}{x_4}$, whereas
$\BR{x_2}{x_3}$ and $\CM{x_1}{x_2}$ come from Lemma~\ref{lem:longconj}\,(ii)
and~(iii). Finally, $\CM{x_2}{x_4}$ is precisely the relation\linebreak
$\cc{b_1}{b_2}{b_3}{b_4}=_G 1$.

The fact that  the two homomorphisms that we have defined are mutually inverse
follows immediately from their definitions.
\end{proof}

\begin{proposition}\label{Prop-RattleEquiv}
For  $1 \leq t \leq n-3$, the quotient $Q_1$ of the Artin group
$A(\Delta_{1,n})$ on generators $b_1,\ldots,b_n$ by the normal
closure of the cycle commutator $\cc{b_1}{b_2}{b_3}{b_4}$ is isomorphic to
the quotient $Q_2$ of the Artin group $A(\Delta_{t,n})$ on generators
	$b_1,\ldots,b_{4+r},c_{5+r},\ldots,c_n$
by the normal closure of the cycle commutator $\cc{b_1}{b_2}{b_3}{b_4}$.
\end{proposition}
\begin{proof} 
	In the homomorphisms that we shall define between these groups,
the generators $b_k$ in the two groups will correspond for $1 \le k \le r+4$,
that is each map will map each such generator $b_k$ of the first group to the generator $b_k$ of the second group,
so there should be no danger of confusion.
In order to find a homomorphism from $Q_2$ to $Q_1$,
we 
	define
elements $c_k$ of $Q_1$  as products of the generators $b_1,\ldots,b_n$ of $Q_1$  follows.
\[c_{5+r} := b_{5+r} b_{4+r} \cdots b_5b_1b_2b_3b_2^{-1}b_1^{-1}b_5^{-1}
 \cdots b_{5+r}^{-1} \]
and $c_k := b_k$ for $6+r \le k \le n$.
Note that, by Lemma~\ref{lem:longconj}\,(i), we have
\[c_{5+r}
=_{Q_1} b_3^{-1}b_2^{-1}b_1^{-1}b_5^{-1}\cdots b_{4+r}^{-1} b_{5+r}b_{4+r}
   \cdots b_5b_1b_2b_3.\]
We only need to verify the relations of $Q_2$ that involve $c_{5+r}$.
Of these $\CM{c_k}{c_{5+r}}$ is immediate for $k > 6+r$.
Since $b_{4+r} \cdots b_5b_1b_2b_3$ commutes with $c_{6+r}=b_{6+r}$
(when $s>1$), $\BR{c_{6+r}}{c_{5+r}}$ follows from $\BR{b_{6+r}}{b_{5+r}}$,
whereas  $\BR{b_3}{c_{5+r}}$ and $\CM{b_k}{c_{5+r}}$ for $k=1,2$ and
$5 \le k \le 4+r$ follow from Lemma~\ref{lem:longconj}\,(ii) and (iii).

Finally, applying Lemma~\ref{lem:ccequiv} to the cycle commutator relation
\linebreak $\cc{b_1}{b_2}{b_3}{b_4} =_{Q_1} 1$
gives $\CM{b_4}{b_1b_2b_3b_2^{-1}b_1^{-1}}$ and,
since $b_4$ commutes with $b_{5+r} \cdots b_5$, we get
$\CM{b_4}{c_{5+r}}$.

In the other direction, in order to define a homomorphism from $Q_1$ to $Q_2$,
we define elements $b_k$ of $Q_2$ for $k \ge 5+r$ as follows.
\[b_{5+r} := b_{4+r} b_{3+r} \cdots b_5b_1b_2b_3c_{5+r}
  b_3^{-1}b_2^{-1}b_1^{-1}b_5^{-1} \cdots b_{4+r}^{-1} \]
and $b_k := c_k$ for $6+r \le k \le n$.
Note that, by Lemma~\ref{lem:longconj}\,(i), we have
\[b_{5+r}
=_{Q_2} c_{5+r}^{-1}b_3^{-1}b_2^{-1}b_1^{-1}b_5^{-1}\cdots b_{3+r}^{-1}
 b_{4+r}b_{3+r} \cdots b_5b_1b_2b_3c_{5+r}.\]
We only need to verify the relations of $Q_1$ that involve $b_{5+r}$.
Of these $\CM{b_k}{b_{5+r}}$ is immediate for $k > 6+r$.
Since $b_{4+r} \cdots b_5b_1b_2b_3$ commutes with $c_{6+r}=b_{6+r}$,
$\BR{b_{6+r}}{b_{5+r}}$ follows from $\BR{c_{6+r}}{c_{5+r}}$, whereas
$\BR{b_{4+r}}{b_{5+r}}$, $\BR{b_{6+r}}{b_{5+r}}$, and
$\CM{b_k}{b_{5+r}}$ for $k=1,2,3$ and
$5 \le k \le 3+r$ follow from Lemma~\ref{lem:longconj}\,(ii) and (iii).

We can also use Lemma~\ref{lem:longconj}\,(i) to give
\[b_{5+r} =_{Q_2} c_{5+r}^{-1} b_{4+r}b_{3+r} \cdots b_5 b_1b_2b_3
b_2^{-1}b_1^{-1}b_5^{-1}\cdots b_{3+r}^{-1} c_{5+r},\]
and then, 
since $c_{5+r}^{-1}b_{4+r}\cdots b_5$ commutes with $b_4$, 
$\CM{b_4}{b_{5+r}}$ follows from $\CM{b_4}{b_1b_2b_3b_2^{-1}b_1^{-1}}$.

It follows immediately from their definitions and Lemma~\ref{lem:longconj}\,(i),
that the two homomorphisms that we have defined are mutually inverse.
\end{proof}

\begin{proofof}{Theorem~\ref{thm:isom}}
Let $a_1,\ldots,a_n$ denote the generators of $G$.
The presentations for $G$ and $Q$ clearly match when $n=4$, so 
we may assume that $n \geq 5$.
	Then we define elements $b_1,\ldots b_{r+4},c_{5+r},c_n$ of $G$ as products of the generators $a_1,\ldots,a_n$, as follows.
\begin{eqnarray*}
b_1 &:=& a_2a_3 \cdots a_{r+1}a_{r+2} a_{r+1}^{-1} \cdots a_3^{-1}a_2^{-1},\\
b_2 &:=& a_2 \cdots a_{r+2} a_{r+3}^{-1} \cdots a_{n-1}^{-1} a_na_{n-1}
 \cdots a_{r+3} a_{r+2}^{-1} \cdots a_2^{-1},\\
 b_3 &:=& a_{r+4}^{-1} \cdots a_{n-1}^{-1} a_na_{n-1} \cdots a_{r+4}\ (=a_n
  \mbox{ when } s=0),\\
b_4 &:=& a_1,\\
b_k &:=& a_{r+7-k} \mbox{ for } 5 \le k \le r+4,\\
c_k &:=& a_{k-1}  \mbox{ for } 5+r \le k \le n\ .
\end{eqnarray*}

We need to verify that the relations of $Q$ hold between the elements
$b_1,\ldots,b_{4+r}$, $c_{r+5},\ldots,c_n$.
The relations of $Q$ among the $b_j$ and $c_k$ that do not involve
$b_1$, $b_2$ or $b_3$ are all immediate, as are $\CM{b_1}{b_3}$,
$\CM{b_1}{c_k}$ for $5+r \le k \le n$
and $\CM{b_3}{b_k}$ for $5 \le k \le r+4$.

Also, since $a_1$ commutes with $a_{n-1} \cdots a_{r+4}$,
$\BR{b_3}{b_4}$ follows immediately from $\BR{a_1}{a_n}$.
Similarly, since
$b_1 =_G a_{r+2}^{-1} \cdots a_3^{-1}a_2a_3 \cdots a_{r+2}$ by
Lemma~\ref{lem:longconj}\,(i), we see that $\BR{b_1}{b_4}$ follows from
$\BR{a_1}{a_2}$ and $\CM{a_1}{a_3 \cdots a_{r+2}}$.

Now $\BR{b_1}{b_5}$, $\CM{b_1}{b_k}$ for $6 \le k \le r+4$, $\BR{b_3}{c_{5+r}}$,
and $\CM{b_3}{c_k}$ for $6+r \le k \le n$ all follow from
Lemma~\ref{lem:longconj}\,(ii) and (iii).

That leaves the relations involving $b_2$.
Let $g := a_{3+r}^{-1} \cdots a_{n-1}^{-1} a_na_{n-1} \cdots a_{r+3}$.
Then, 
since \[ b_1=a_2\cdots a_{r+1}a_{r+2}a_{r+1}^{-1}\cdots a_2^{-1}
\ \mbox{and}\  
b_2=a_2\cdots a_{r+1}a_{r+2}ga_{r+2}^{-1}a_{r+1}^{-1}\cdots a_2^{-1},\] 
$\BR{b_1}{b_2}$ is equivalent to $\BR{a_{r+2}}{a_{r+2}ga_{r+2}^{-1}}$
and hence to $\BR{a_{r+2}}{g}$, which follows from
Lemma~\ref{lem:longconj}\,(ii).

To verify $\BR{b_2}{b_3}$, again note
that $b_2 = a_2\cdots a_{r+2} g a_{r+2}^{-1} \cdots a_2^{-1}$
Then, since $b_3 =_G a_{r+3}ga_{r+3}^{-1}$ and $a_2, \ldots, a_{r+2}$
commute with $b_3$, we see that $\BR{b_2}{b_3}$ is equivalent to
$\BR{g}{a_{r+3}ga_{r+3}^{-1}}$.
Now Lemma~\ref{lem:longconj}\,(ii) gives $\BR{a_{r+3}}{g}$, so
$a_{r+3}ga_{r+3}^{-1} =_G g^{-1}a_{r+3}g$,
from which $\BR{g}{a_{r+3}ga_{r+3}^{-1}}$
(and hence also $\BR{b_2}{b_3}$) follows immediately. 

By Lemma~\ref{lem:tcequiv}\,(ii), the relation $\CM{b_2}{b_4}$ follows from
the relator\linebreak
$\tc{a_1}{a_2}{\ldots}{a_n}_{s+1}$ of $G$.

As we saw above, Lemma~\ref{lem:longconj}\,(ii) gives
$\BR{a_{r+2}}{g}$ and, since we also have $\CM{a_i}{g}$
for $2 \le i \le r+1$, we can apply Lemma~\ref{lem:longconj} to the sequence
$a_2,\ldots,a_{r+2},g$, and then Lemma~\ref{lem:longconj}\,(iii) implies
$\CM{b_2}{b_k}$ for $5 \le k \le r+4$.
Furthermore 
(since $a_2\cdots a_{r+2}$ commutes with $c_k$),
$\CM{b_2}{c_k}$ is equivalent to $\CM{g}{c_k}$
for $5+r \le k \le n$, which also follows from Lemma~\ref{lem:longconj}\,(iii).

It remains to verify that $\tc{b_1}{b_2}{b_3}{b_4}_1 =_G 1$.
Since $b_3$ commutes with $a_k$ for $2 \le k \le a_{r+2}$, and
$\BR{b_3}{a_{r+3}}$ by Lemma~\ref{lem:longconj}\,(ii), while also
$b_2=a_2 \cdots a_{r+2}(a_{r+3}^{-1}b_3a_{r+3})a_{r+2}^{-1} \cdots a_2^{-1}$,
we have
\begin{eqnarray*}
b_2b_3b_2^{-1} &=_G&  a_2 \cdots a_{r+2} (a_{r+3}^{-1} b_3 a_{r+3})
b_3 (a_{r+3}^{-1} b_3^{-1} a_{r+3}) a_{r+2}^{-1} \cdots a_2^{-1}\\ 
&=_G& a_2 \cdots a_{r+2} a_{r+3} a_{r+2}^{-1} \cdots a_2^{-1} 
\end{eqnarray*}
Now, using this expression for $b_2b_3b_2^{-1}$, we find that
$b_1^{-1}b_2b_3b_2^{-1}b_1$ reduces to $a_{r+3}$
using free reduction and the commuting relations
between $a_{r+3}$ and $a_2,\ldots,a_{r+1}$.
Hence
$\CM{b_4}{b_1^{-1}b_2b_3b_2^{-1}b_1}$ follows from $\CM{a_1}{a_{r+3}}$,
and now $\tc{b_1}{b_2}{b_3}{b_4}_1 =_G 1$ follows from this,
by Lemma~\ref{lem:tcequiv}\,(i).

In the other direction, let $b_1,\ldots,c_n$ denote the generators of $Q$,
and define
elements $a_1,\ldots,a_n$ in $Q$ as follows.

\begin{eqnarray*}
&&a_1:=b_4,\\
&&a_2 := b_{r+4}b_{r+3} \cdots b_5b_1b_5^{-1} \cdots b_{r+3}^{-1}b_{r+4}^{-1},\\
&&a_k := b_{r+7-k}  \mbox{ for } 3 \le k \le r+2,\\
&&a_{r+3} := b_1^{-1}b_2b_3b_2^{-1}b_1,\\
&& a_k := c_{k+1}  \mbox{ for } r+4 \le k \le n-1,\\
&& a_n := c_n \cdots c_{r+5} b_3 c_{r+5}^{-1} \cdots c_n^{-1}\ (=b_3
  \mbox{ when } s=0).
\end{eqnarray*}

We need to verify that the relations of $G$ hold between the elements $a_1,\ldots,a_n$ of $Q$.

The relations of $G$ among the $a_i$ that do not involve
$a_2$, $a_{r+3}$ or $a_n$
are all immediate, as are $\CM{a_2}{a_k}$ for $r+4 \le k \le n-1$,
$\CM{a_n}{a_k}$ for $3 \le k \le r+2$,
and $\CM{a_{r+3}}{a_k}$ for $3 \le k \le r+1$ and $r+5 \le k \le n-1$.

Now $\BR{a_1}{a_n}$ follows from $\BR{b_4}{b_3}$
since $b_4$ commutes with $c_n \cdots c_{r+5}$, and
$\CM{a_2}{a_n}$ follows similarly from $\CM{b_1}{b_3}$, whereas
$\CM{a_1}{a_{r+3}}$ comes from $\tc{b_1}{b_2}{b_3}{b_4}=_G 1$
and Lemma~\ref{lem:tcequiv}\,(i).

Since $a_{r+3}$ commutes with $b_{r+4} \cdots b_6$,
$\CM{a_2}{a_{r+3}}$ reduces to $\CM{b_5b_1b_5^{-1}}{a_{r+3}}$, which is
equivalent to $\CM{b_1^{-1}b_5b_1}{a_{r+3}}$, and this follows from
$\CM{b_5}{b_2b_3b_2^{-1}}$.
Similarly, since $a_1$ commutes with $b_{r+4} \cdots b_5$, 
$\BR{a_1}{a_2}$ reduces to $\BR{b_4}{b_1}$.

Next we see that
$\BR{a_2}{a_3}$ comes from Lemma~\ref{lem:longconj}\,(ii), and $\CM{a_2}{a_k}$
for $4 \le k \le r+2$ come from  Lemma~\ref{lem:longconj}\,(iii).
Similarly, $\BR{a_{n-1}}{a_n}$ comes from Lemma~\ref{lem:longconj}\,(ii) and
$\CM{a_k}{a_n}$ for $r+4 \le k \le n-2$ come from
Lemma~\ref{lem:longconj}\,(iii).

Furthermore, we have $\BR{b_1}{b_2b_3b_2^{-1}}$ by Lemma~\ref{lem:brcm}\,(i)
and, since we also have $\CM{b_5}{b_2b_3b_2^{-1}}$, we get
$\BR{a_{r+2}}{a_{r+3}} = \BR{b_5}{b_1^{-1}b_2b_3b_2^{-1}b_1}$ by applying
Lemma~\ref{lem:brcm}\,(i) with $f=b_5$, $g=b_1$ and $h=b_2b_3b_2^{-1}$.


From $\BR{c_{r+5}}{b_3}$  and $\CM{c_{r+5}}{b_1^{-1}b_2}$, we get
$\BR{c_{r+5}}{b_1^{-1}b_2b_3b_2^{-1}b_1} = \BR{a_{r+4}}{a_{r+3}}$.

Since $b_1$ commutes with $a_n$ and $b_2b_3b_2^{-1}$ commutes with
$c_n \cdots c_{r+6}$, $\CM{a_n}{a_{r+3}}$ reduces to
$\CM{c_{r+5} b_3 c_{r+5}^{-1}}{b_2b_3b_2^{-1}}$, which we get by applying
Lemma~\ref{lem:brcm}\,(ii) with $f=c_{r+5}$, $g=b_3$ and $h=b_2$.

It remains to check that $\tc{a_1}{a_2}{\ldots}{a_n}_{s+1} =_G 1$.
We find that
$$a_2 \cdots a_{r+2} a_{r+3}^{-1} \cdots a_{n-1}^{-1} a_n a_{n-1}
 \cdots a_{r+3} a_{r+2}^{-1}\cdots a_2^{-1}$$
freely reduces to $b_{r+4} \cdots b_5 b_2b_3^{-1}b_2^{-1} b_1 b_3 b_1^{-1}
b_2 b_3 b_2^{-1} b_5^{-1} \cdots b_{r+4}^{-1}$ which,
since $b_1$ and $b_3$ commute, is equal in $G$ to
\begin{eqnarray*}
b_{r+4} \cdots b_5 b_2b_3^{-1}b_2^{-1} b_3 b_2 b_3 b_2^{-1} b_5^{-1}
	\cdots b_{r+4}^{-1} &=_G&\\
b_{r+4} \cdots b_5 b_2 b_5^{-1} \cdots b_{r+4}^{-1} &=_G & b_2,
\end{eqnarray*}
and this commutes with $a_1 = b_4$, so the relation
$\tc{a_1}{a_2}{\ldots}{a_n}_{s+1} =_G 1$ follows from
Lemma~\ref{lem:tcequiv}\,(ii).

Finally, we need to check that both composites of the homomorphisms between
the two groups are equal to the identity.

First consider the images of the $a_i$ under the composite mapping
$G \to Q \to G$. They map back to $a_i$ immediately except when
$i=2,r+3$ or $n$, and for $a_n$ this follows using a free reduction.
Note that by Lemma~\ref{lem:longconj}\,(i), the definition of
$b_1$ in $G$ is equivalent to
$b_1 := a_{r+2}^{-1} \cdots a_3^{-1}a_2a_3 \cdots a_{r+2}$ and using this,
we find that $a_2$ also maps back to itself by using a free reduction.
As for $a_{r+3}$, we saw earlier that the image of $b_1^{-1}b_2b_3b_2^{-1}b_1$
under the map $Q \to G$ reduces to $a_{r+3}$ in $G$ by using free reduction and the commuting relations, and so it too maps back to itself.

The images of $b_i$ under the composite mapping $Q \to G \to Q$ are immediately
equal to $b_i$ except for $b_1$, $b_2$ and $b_3$. The same applies to $b_3$
using a free reduction, and also for $b_1$ after rewriting it as
$a_{r+2}^{-1} \cdots a_3^{-1}a_2a_3 \cdots a_{r+2}$ as above.

Finally, for $b_2$, note first that the image of $a_2 \cdots a_{r+2}$ under
the map $G \to Q$ freely reduces to $b_{r+4} \cdots b_5b_1$,
whereas the image of $a_{r+4}^{-1} \cdots a_{n-1}^{-1} a_na_{n-1}
 \cdots a_{r+4}$ freely reduces to $b_3$, and hence that of
$a_{r+3}^{-1} \cdots a_{n-1}^{-1} a_na_{n-1} \cdots a_{r+3}$ reduces to 
\[
b_1^{-1}b_2b_3^{-1}b_2^{-1}b_1b_3b_1^{-1}b_2b_3b_2^{-1}b_1
=_{Q}  b_1^{-1}b_2b_3^{-1}b_2^{-1}b_3b_2b_3b_2^{-1}b_1
=_{Q} b_1^{-1} b_2 b_1
\]
and, since $b_2$ commutes with $b_{r+4}\cdots b_5$, we find that the image
of $b_2$ under the composite is indeed $b_2$.
\end{proofof}

\section{Proving non-isomorphism}\label{sec:nonisoms}
The aim of this section is to prove Theorem~\ref{thm:nonisom}.
By Proposition~\ref{PropArtin-n-gon} and Theorem~\ref{thm:isom}, this is
equivalent to proving that the quotients of $A(\Delta_n)$ by the normal
closures of the cycle and twisted cycle commutators
$\cc{a_1}{a_2}{\ldots}{a_n}$ and 
$\tc{a_1}{a_2}{\ldots}{a_n}_t$ are not isomorphic for any $t$ with
$1 \leq t \leq n-3$, which is what we shall prove here.
(Note that these quotients are isomorphic  when $t=n-2$, which can be
seen by considering the automorphism of $A(\Delta_n)$ induced by
$a_i \mapsto a_i^{-1}$ for $1 \leq i \leq n$.)

When $n=4$, the result is easily proved by computer. For example, the two
quotients have $9$ and $8$ conjugacy classes of subgroups of index $4$.
So we shall assume from now on that $n > 4$.

We start with a result that is probably already known. 
In general, two surjective homomorphisms $\sigma_1,\sigma_2: G \to H$ are said
to be \emph{equivalent} if they have the same kernels or, equivalently, if
there is an automorphism $\alpha$ of $H$ with
$\alpha(\sigma_1(g)) = \sigma_2(g)$ for all $g \in G$.

\begin{proposition}\label{prop:antosn}
All surjective homomorphisms $\sigma:A(\Delta_n) \to \Sym(n)$ are equivalent
when $n > 4$.
\end{proposition}
\begin{proof} Let $\sigma:A(\Delta_n) \to \Sym(n)$ be a surjective homomorphism.
We consider the restriction of $\sigma$ to the subgroup
$B := \langle a_1,a_2,\ldots,a_{n-1} \rangle$ of $A(\Delta_n)$, which is
isomorphic to the braid group $B_n$. If $\sigma(B)$ is abelian then, since
$a_2$ and $a_n$ commute, we have $\sigma(a_2) \in Z(\Sym(n)) = 1$. But the
generators $a_i$ are all conjugate in $A(\Delta_n)$ so this is impossible.
Hence $\sigma(B)$ is nonabelian.

In~\cite[Theorem A]{Lin} Lin proved that, for $n >4$, all homomorphisms of
$B_n$ to $\Sym(k)$ with $k < n$ have cyclic image. So, if $\sigma(B)$ is
intransitive then it is abelian, contrary to what we just proved.
Artin proved in \cite{Artin} that, for $n > 4$, all homomorphisms from $B_n$
to a transitive subgroup of $\Sym(n)$ are surjective and equivalent to
the homomorphism in which the generators map to the transpositions
$(i,i+1)$ for $1 \leq i < n$.

So we may assume that $\sigma(a_i) = (i,i+1)$ for $1 \leq i < n$.
Now, since $\sigma(a_n)$ commutes with the images of $\sigma(a_i)$ for
$2 \leq i \leq n-2$ and these images generate the subgroup $\Sym(n-2)$ of
$\Sym(n)$ acting on $\{2,3,\ldots,n-1\}$, which has trivial centre,
we see that the only possible image of $\sigma_n$ is $(1,n)$.
\end{proof}

\subsection{Subgroups of $A(\Delta_n)$ and its quotients}\label{subsec:hi}
Let $G :=  A(\Delta_n) = \langle X \mid R \rangle$ be our
standard presentation of $\Delta_n$ with $X = \{ a_1,\ldots,a_n \}$.
Let $G_0$ and $G_t$ for $1 \leq t \leq n-2$ be the quotients of $G$ by the normal
closures of the cycle commutator
$\cc{a_1}{a_2}{\ldots}{a_n}$ and the twisted cycle commutator
$\tc{a_1}{a_2}{\ldots}{a_n}_t$, respectively.
Define $\sigma:G \to \Sym(n)$ by $\sigma(a_i) = (i,i+1)$ for $1 \leq i < n$
and $\sigma(a_n) = (1,n)$. Then, since we are assuming that $n>4$,
Proposition~\ref{prop:antosn} tells us that $\sigma$ is a representative
of the unique equivalence class of surjective homomorphisms $G \to \Sym(n)$. 

Since $\cc{a_1}{a_2}{\ldots}{a_n}$ and $\tc{a_1}{a_2}{\ldots}{a_n}_t$
are both contained in $\ker(\sigma)$, it follows that there is also a unique
equivalence class of surjective homomorphisms $G_i \to \Sym(n)$ for
$0 \leq i \leq n$, with the map $\sigma_i$ induced by $\sigma$ as representative.

Let $H < G$ be the inverse image in $G$ of the stabiliser in $\Sym(n)$ of the
unordered pair $\{1,2\}$ under $\sigma$.  So $|G:H| = n(n-1)/2$.
Define $H_i$ for $1 \leq i \leq n-2$ to be the corresponding subgroups of $G_i$,
and note that they also have index $n(n-1)/2$ in $G_i$. In fact the
uniqueness of the equivalence classes of surjective homomorphisms
$G_i \to \Sym(n)$ implies that an isomorphism $G_0 \to G_t$ for
$1 \leq t \leq n-2$
would map $H_0$ to an image of $H_t$ under an automorphism of $\Sym(n)$, and so
there would be such an isomorphism mapping $H_0$ to $H_t$.  The remainder of
this section will be devoted to the proof of Proposition~\ref{prop:aqhi}, stated below, which shows
that $H_0$ and $H_t$ have different abelianisations when $1 \leq t \leq n-3$,
and so they cannot be isomorphic.  This will also complete the proof of
Theorem~\ref{thm:nonisom}.

\begin{remark} We initially carried out corresponding calculations using
the inverse image in $G$ of the stabiliser of $1$ in $\Sym(n)$, which has
index $n$ in $G$, but we found that the corresponding subgroups $H_0$ and $H_t$
both had the abelianisation $\frac{\Z}{2\Z}\oplus \Z^2$, so that did not work.
\end{remark} 

\begin{proposition}\label{prop:aqhi}
For all $n > 4$
we have $H \cong \Z^4$, $H_0/[H_0,H_0] \cong \frac{\Z}{2\Z}\oplus \Z^3$, and
$H_t/[H_t,H_t] \cong \frac{\Z}{2\Z} \oplus \frac{\Z}{4\Z} \oplus \Z^2$
for $1 \leq t \leq n-3$.
\end{proposition}
The proof of this proposition will be completed only at the very end of 
Section~\ref{sec:tabHrels}.

\subsection{Computing presentations of subgroups} We shall now present a quick
summary of the methods based on the Reidemeister-Schreier algorithm that
are used in computer calculations of presentations of subgroups of finite index
of groups defined by a finite presentation.
In the following subsections we shall use this theory to carry out a calculation
of this type by hand for an infinite family of examples, although this
calculation was based on observations of the results of computer calculations
in small cases.  This theory can be found in many textbooks; our description
here is based on that in \cite[Section 2.5.2]{HEO}. 

Let $\GG := \langle X \mid R \rangle$ be a group defined by a presentation,
let $\HH \leq \GG$, and let $T$ be a right transversal of $\HH$ in $\GG$ that contains
the empty word $\emptyword$. For $t \in T$ and $g \in \GG$, we denote the unique
element of $\HH tg \cap T$ by $\overline{tg}$. (So $\overline{tg}$ is the image of
$t$ under $g$ in the action of $\GG$ on $T$ induced by its natural action on the
right cosets of $\HH$ in $\GG$.)
In the remainder of this section, we shall apply the results of this section to
various pairs $\GG,\HH$ that have just been introduced, namely to the pair
$(\GG,\HH)=(G,H)=(A(\Delta_n),H)$ as well as the pairs $(G_i,H_i)$ of
subgroups of those.  Hence we have been careful to distinguish between
$(\GG,\HH)$ and $(G,H)$ typographically.

Let $Y$ be a subset of $\HH$ and suppose that, for each $t \in T$ and $x \in X$,
there exists a word $\rho(t,x) \in (Y^\pm)^*$ with $tx =_\GG \rho(t,x) u$, where
$u = \overline{tx}$. This implies that $u x^{-1} =_\GG \rho(t,x)^{-1} t$,
and we define $\rho(u,x^{-1}) := \rho(t,x)^{-1}$.

We can now recursively extend the definition of $\rho$ to words
$w \in (X^\pm)^*$, by defining $\rho(t,\emptyword) = \emptyword$
for all $t \in T$ and, for a word $v \in (X^\pm)^*$ and $x \in X^{\pm}$,
$\rho(t,vx) = \rho(t,v)\rho(u,x)$, where $\overline{tv} = u$. So we have
$tw =_\GG \rho(t,w) \overline{tw}$ for all $t \in T$ and $w \in (X^\pm)^*$. 

Note that, if $w$ represents an element of $\HH$, then $\overline{w}= \emptyword$,
and so $w =_\GG \rho(t,w)$, which proves:

\begin{proposition}\label{prop:Hgen} Under the assumptions above,
we have $\HH = \langle Y \rangle$.
\end{proposition}

Suppose that, for each $y \in Y$, we are given a word
$\varphi(y) \in (X^{\pm})^*$ with $y =_\GG \varphi(y)$.
The following result, which we shall use in our hand calculations, is proved
in \cite[Theorem 2.63]{HEO} and the remark that follows it.
Notice that, we are following the common practice of abusing notation by
using $Y$ to denote both a subset of $\HH$ and a set of generators in a
presentation of $\HH$.

\begin{theorem}\label{thm:subpres}
Under the assumptions above, $\langle Y \mid S_1 \cup S_2 \rangle$ is a
presentation of $\HH$ on the generating set $Y$, where
$S_1 = \{\rho(t,w) : t \in T, w \in R\}$, and
$S_2 = \{ \rho(\epsilon,\varphi(y))y^{-1} : y \in Y \}$.
\end{theorem}

In our application below, the relators $R$ are given as a set $R'$ of
relations $w_1 = w_2$,
and it is convenient to replace the relators in $S_1$ by the equivalent
relations $\{\rho(t,w_1) = \rho(t,w_2) : t \in T, (w_1 = w_2)  \in R' \}$.

\subsection{The words $\rho(t,x)$ for the subgroup $H$ of $G$.}
Our aim in the remainder of this section is to apply Theorem~\ref{thm:subpres}
to the subgroups $H$ of $G$ and $H_i$ of $G_i$ that were defined in Subsection~\ref{subsec:hi}.
We start by finding a generating set $Y$ of $H$ and a transversal $T$ of $H$ in
$G$, and computing the words $\rho(t,a_i) \in (Y^\pm)^*$ for $t \in T$ and
$1 \leq i \leq n$.

Define the words
\[ \xi_1 := a_1,\  \xi_i := a_{i+1}\ (2 \leq i \leq n-2),\ 
   \xi_{n-1} := a_2^2,\ \xi_n := a_n^2,\ \xi_{n+1} := a_2a_1\prod_{i=3}^n a_i.\]

It is straightforward to check that the words $\xi_i \in (X^\pm)^*$ represent
elements of $H$ for $1 \leq i \leq n+1$, and we denote the element of $H$ that is
represented by the word $\xi_i$ by $y_i$.
So, in the notation of the preceding subsection, we can define
$\varphi(y_i) := \xi_i$. It follows from Proposition~\ref{prop:ect} below and
Theorem~\ref{thm:subpres} that $H = \langle y_i: 1 \leq i \leq n+1  \rangle$.

For $1 \leq k < l \leq n$, define
\[ t_{k,l} := \prod_{i=2}^{l-1} a_i \prod_{i=1}^{k-1} a_i. \]
Then the element $\sigma(t_{k,l}) \in \Sym(n)$ maps $1$ to $k$ and
$2$ to $l$ (note that we are composing permutations from left to right),
so the elements $t_{k,l}$ form a right transversal of $H$ in $G$.

\begin{proposition}\label{prop:ect}
For $1 \leq k < l \leq n$ and $ 1 \leq m \leq n$, 
let words $\rho(t_{k,l},a_m) \in Y^*$ and pairs of integers $k',l'$ be as specified by the following table.

Then for all such $k,l,m$, the equation
$t_{k,l}a_m =_G \rho(t_{k,l},a_m) t_{k',l'}$ holds.
\begin{center}

\begin{tabular}{|l|ll|}
\hline
	Case & $\rho(t_{k,l},a_m)$ & $(k',l')$ \\
\hline
$l<m < n$ & $y_{m-1}$ & $(k,l)$\\
$l=m < n$ & $\emptyword$ & $(k,l+1)$\\
$k=m =l-1$ & $y_1$ & $(k,l)$\\
$k=m<l-1$ & $\emptyword$ & $(k+1,l)$\\
$k<m =l-1$ & $y_{n-1}^{\prod_{i=2}^{l-2} y_i}$ & $(k,l-1)$\\
$k<m <l-1$ & $y_m$ & $(k,l)$\\
$k =m+1$ & $y_{n-1}^{\prod_{i=1}^{k-1} y_i}$ & $(k-1,l)$\\
$k > m+1$ & $y_{m+1}$ & $(k,l)$\\
$k=1, l=m=n$ & $y_1^{y_{n+1}^{-1}}$ & $(k,l)$\\
$k=1, l<m=n$ & $y_n y_{n+1}^{-1}$ & $(l,n)$\\
$k>1, l=m=n$ & $y_{n+1}$ & $(1,k)$\\
$k>1, l<m=n$ & $y_{n-2}^{y_{n+1}^{-1}}$ & $(k,l)$\\
\hline
\end{tabular}
\end{center}
\end{proposition}

\begin{proof} We shall go through the table entries line by line, proving the
claimed result in each case.

When $l<m<n$, $a_m$ commutes with $t_{k,l}$ and $a_m = y_{m-1}$.

When $l=m<n$, $a_m$ commutes with $\prod_{j=1}^{k-1} a_i$, and so
$t_{k',l'} =_G t_{k,l+1}$.

When $k=m=l-1$, we need to prove that
\[\prod_{i=2}^{l-1} a_i  \left(\prod_{j=1}^{l-2} a_j\right)a_{l-1} =_G
a_1 \prod_{i=2}^{l-1} a_i \prod_{j=1}^{l-2} a_j, \]
which we do by induction on $l$. The base case $l=2$ is trivial.
Since $a_{l-1}$ commutes with $a_1 \cdots a_{l-3}$ we have, using induction,
\begin{eqnarray*}
 \prod_{i=2}^{l-1} a_i \left(\prod_{j=1}^{l-2} a_j\right)a_{l-1} &=_G&
  \prod_{i=2}^{l-2} a_i \left(\prod_{j=1}^{l-3} a_j\right)
    a_{l-1}a_{l-2}a_{l-1} =_G\\
\prod_{i=2}^{l-2} a_i \left(\prod_{j=1}^{l-3} a_j\right)a_{l-2}a_{l-1}a_{l-2}
    &=_G&
a_1 \prod_{i=2}^{l-2} a_i \left(\prod_{j=1}^{l-3} a_j\right) a_{l-1}a_{l-2}
    =_G\\
&&a_1 \prod_{i=2}^{l-1} a_i \prod_{j=1}^{l-2} a_j.
\end{eqnarray*}

The result for $k=m < l-1$ is immediate.

	When $k < m = l-1$,  note first that the claimed value of $\rho(t_{k,l},a_m)$ is
\[ y_{n-1}^{\prod_{i=2}^{l-2} y_i} = (a_2^2)^{\prod_{i=3}^{l-1} a_i} =_G
(a_2^{\prod_{i=3}^{l-1} a_i})^2 \]
which, by Lemma~\ref{lem:longconj}, is equal in $G$ to
$(a_{l-1}^2)^{ \left({\prod_{i=2}^{l-2} a_i}\right)^{\!-1}}$ and so,
after free cancellation, we have
$y_{n-1}^{\prod_{i=2}^{l-2} y_i} t_{k,l-1} =_G
   \left(\prod_{i=2}^{l-2} a_i \right) a_{l-1}^2 \prod_{j=1}^{k-1} a_j$,
and the claim follows from the fact that $a_{l-1}$ commutes with
$\prod_{j=1}^{k-1}a_j$. 

When $k < m < l-1$, since $a_m$ commutes with $\prod_{j=1}^{k-1} a_j$ and
with $\prod_{i=m+2}^{l-1} a_i$, the proof reduces to showing that
$\left(\prod_{i=2}^{m+1}a_i\right)a_m = a_{m+1} \prod_{i=2}^{m+1}a_i$,
which follows from
$\BR{a_m}{a_{m+1}}$ and the fact that $a_{m+1}$ commutes with
$\prod_{i=2}^{m-1}a_i$.  

The case $k=m+1$ seems to be the most difficult to prove. Note first that
$a_1$ commuting with $\prod_{i=3}^k a_i$ followed by Lemma~\ref{lem:longconj}
give
\[y_{n-1}^{\prod_{i=1}^{k-1} y_i} =_G
a_1^{-1} (a_2^2)^{\left({\prod_{i=3}^k a_i}\right)} a_1 =_G
a_1^{-1} (a_k^2)^{\left({\prod_{i=2}^{k-1} a_i}\right)^{\!-1}} a_1,\] so
\[y_{n-1}^{\prod_{i=1}^{k-1} y_i} t_{k-1,l} =_G
a_1^{-1} (a_k^2)^{{\left(\prod_{i=2}^{k-1} a_i\right)}^{\!-1}} a_1
  \prod_{i=2}^{l-1}a_i \prod_{j=1}^{k-2} a_j.\]
Now, the right hand side of this equality contains the subword
$a_1^{\prod_{i=2}^{k-1} a_i}$ which, by Lemma~\ref{lem:longconj} is equal in
$G$ to $a_{k-1}^{\left( \prod_{i=1}^{k-2} a_i\right)^{\!-1}}$, and this 
substitution results in the expression
\[ a_1^{-1} \left( \prod_{i=2}^{k-1} a_i \right) a_k^2\,
a_{k-1}^{\left( \prod_{i=1}^{k-2} a_i\right)^{\!-1}}
 \prod_{i=k}^{l-1}a_i \prod_{j=1}^{k-2} a_j\]
which, since $\prod_{i=1}^{k-2} a_i$ commutes with $\prod_{i=k}^{l-1}a_i$,
reduces to
\( a_1^{-1} \left( \prod_{i=2}^{k-1} a_i \right) a_k^2 
 \prod_{i=1}^{l-1}a_i.\)
Now, by applying Lemma~\ref{lem:hardcase} below to a prefix of this expression
followed by commutativity relations, we find that it is equal in $G$ to
\[ \prod_{i=2}^{k-1}a_i \left(\prod_{i=1}^{k-2}a_i\right) a_ka_{k-1}^2
  \prod_{i=k+1}^{l-1}a_i =_G
\prod_{i=2}^{l-1}a_i \left(\prod_{j=1}^{k-2}a_j\right) a_{k-1}^2 =
  t_{k,l} a_{k-1}=t_{k,l}a_m.\]

When $k>m+1$ we find, by using the commutativity relations and
$\BR{a_m}{a_{m+1}}$, that $\left(\prod_{j=1}^{k-1}a_j\right) a_m =_G
a_{m+1} \prod_{j=1}^{k-1}a_j$ and similarly, since $m+1 < l-1$, we have
$\left(\prod_{j=2}^{l-1}a_j\right) a_{m+1} =_G a_{m+2} \prod_{j=2}^{l-1}a_j$.
So in this case \[ t_{k,l}a_m =_G a_{m+2}t_{k,l}=y_{m+1}t_{k,l} \]
and the result follows.

For the case $k=1$ and $l=m=n$, note that $y_1^{y_{n+1}^{-1}} =
a_1^{{\left(a_2a_1\prod_{i=3}^na_i\right)}^{\!-1}}$ and, since $a_1$
commutes with $\prod_{i=3}^{n-1}a_i$, this is equal in $G$ to
\[(a_1a_na_1a_n^{-1}a_1^{-1})^{{\left(\prod_{i=2}^{n-1}a_i\right)}^{\!-1}} =_G
a_n^{{\left(\prod_{i=2}^{n-1}a_i\right)}^{\!-1}}\] (using $\BR{a_1}{a_n}).$
It follows by free cancellation that 
\[y_1^{y_{n+1}^{-1}}t_{1,n} = y_1^{y_{n+1}^{-1}}\prod_{i=2}^{n-1}a_i 
= \left( \prod_{i=2}^{n-1}a_i\right) a_n=t_{1,n}a_m.\]

When $k=1$ and $l<m = n$, we apply commutativity to get
\[y_ny_{n+1}^{-1}t_{l,n} =
a_n^2\left(\prod_{i=3}^{n}a_i\right)^{\!\!\!-1}
a_1^{-1}a_2^{-1} \prod_{i=2}^{n-1}a_i \prod_{j=1}^{l-1}a_j
=_G a_na_1^{-1} \left(\prod_{i=2}^{n-1}a_i\right)^{\!\!\!-1}
\prod_{i=2}^{n-1}a_i \prod_{j=1}^{l-1}a_j.\]
The final product freely reduces to $a_n  \prod_{i=2}^{l-1}a_i$ which, since $a_n$
commutes with $\prod_{i=2}^{l-1}a_i$, is equal in $G$ to $t_{1,l}a_n$.

For the case $k>1$, $l=m=n$, we use first that $a_n$ commutes with
commutes with $\prod_{j=2}^{k-1} a_j$ and then
that $a_1$ commutes
with $\prod_{i=3}^{n-1} a_i$ 
to see that
\begin{eqnarray*}
	t_{k,l}a_m&=&t_{k,n}a_n=\prod_{i=2}^{n-1}a_i\left(\prod_{i=1}^{k-1}a_i\right) a_n 
	=_G\left( \prod_{i=2}^{n-1}a_1\right) a_n\left(\prod_{i=2}^{k-1}a_i\right)\\
	&=& \left(a_2a_1 \prod_{i=3}^na_i\right ) \prod_{i=2}^{k-1}a_i = y_{n+1}t_{1,k}.
\end{eqnarray*}

Finally, for the case $k>1$, $l < m=n$, 
we need to verify
\[ 
	\left( \prod_{i=2}^{l-1}a_i \prod_{i=1}^{k-1} a_i\right)a_n = _G
	y_{n-2}^{y_{n+1}^{-1}}\left( \prod_{i=2}^{l-1}a_i \prod_{i=1}^{k-1} a_i\right) \]
We note that, since $a_1$ commutes
with $\prod_{i=3}^{n-2} a_i$,
we have \begin{eqnarray*}
	y_{n-2}^{y_{n+1}^{-1}} &=_G&
a_{n-1}^{\left ( a_2a_1 \prod_{i=3}^na_i\right )^{-1}}
	=_G
a_{n-1}^{\left (  \left ( \prod_{i=2}^{n-2}a_i\right ) a_{n-1}a_1a_n\right )^{-1}}\\
	&=_G&
\left( a_{n-1}^{\left (a_{n-1}a_1a_n\right)^{-1}}\right)^{ \left ( \prod_{i=2}^{n-2}a_i\right )^{-1}}
= \left (a_{n-1}a_1a_na_{n-1}a_n^{-1}a_1^{-1}a_{n-1}^{-1}\right)
	       ^{ \left ( \prod_{i=2}^{n-2}a_i\right )^{-1}}
\end{eqnarray*}
Braid and commutator relations reduce that last product to 
\[ \left (a_1a_na_1^{-1}\right)
	       ^{ \left ( \prod_{i=2}^{n-2}a_i\right )^{-1}},\]
	       and then the commutator relators between $a_i$ for $i=l,\ldots,n-2$ and both $a_1$ and $a_n$ reduce it further to 
\[ (a_1a_na_1^{-1})
	       ^{ \left ( \prod_{i=2}^{l-1}a_i\right )^{-1}}.\]
	       So now, applying commutator relations, we deduce that
\begin{eqnarray*}
	y_{n-2}^{y_{n+1}^{-1}}t_{k,l} &=_G& 
\left ( \prod_{i=2}^{l-1}a_i\right )\left (a_1a_na_1^{-1}\right)
	                          \left ( \prod_{i=1}^{k-1}a_i\right )
=_G \left ( \prod_{i=2}^{l-1}a_i\right )\left (a_1a_n\right)
	                          \left ( \prod_{i=2}^{k-1}a_i\right )\\
				  &=_G&
		       \left ( \prod_{i=2}^{l-1}a_i 
		        \prod_{i=1}^{k-1}a_i\right ) a_n 
=_G t_{k,l}a_n \end{eqnarray*}
\end{proof}

\begin{lemma}\label{lem:hardcase}
Suppose that $1 < k < n$. Then
\[ a_1^{-1}\left(\prod_{i=2}^{k-1}a_i\right) a_k^2 \prod_{i=1}^{k} a_i =_G
\prod_{i=2}^{k-1}a_i \left(\prod_{i=1}^{k-2}a_i\right) a_ka_{k-1}^2.\]
\end{lemma}
\begin{proof} The proof is by induction on $k$. In the base case $k=2$, the
relation to be proved is
$a_1^{-1}a_2^2a_1a_2 =_G a_2a_1^2$, which follows from $\BR{a_1}{a_2}$.
So assume that $k>2$.

Using the fact that $a_1$ commutes with $a_k^2$ and with $\prod_{i=3}^{k-1}a_i$,
as well as $\BR{a_1}{a_2}$, we deduce that the left hand side of the required relation is equal in
$G$ to
\[ a_1^{-1}a_2a_1 \left(\prod_{i=3}^{k-1}a_i\right) a_k^2 \prod_{i=2}^{k} a_i
=_G a_2a_1a_2^{-1} \left(\prod_{i=3}^{k-1}a_i\right) a_k^2 \prod_{i=2}^{k} a_i
\]
which, by applying the inductive hypothesis applied to the sequence
$a_2,a_3,\ldots,a_k$ followed by commutativity relations, is equal in $G$ first to
	\[a_2a_1 \prod_{i=3}^{k-1}a_i \left(\prod_{i=2}^{k-2}a_i\right) a_ka_{k-1}^2\]
	and then (since $a_1$ commutes with each $a_i$ with $i\geq 3$) to
	\[\prod_{i=2}^{k-1}a_i \left(\prod_{i=1}^{k-2}a_i\right) a_ka_{k-1}^2.\]
\end{proof}

\subsection{The relations of $H/[H,H]$.}
\label{sec:abHrels}
We shall now use Theorem~\ref{thm:subpres} to calculate the set
$S_1 \cup S_2$ of defining relators of
$H$ on the generators $y_i$, using in those calculations the values of
$\rho(t_{k,l},a_m)$ that were verified in Proposition~\ref{prop:ect}.
Since we want to calculate a presentation of
$H/[H,H]$, we shall abelianise the relators of $H$ immediately, and write them as words in
the images $z_i$ of $y_i$ in $H/[H,H]$ using additive notation. Note that this
means that each of the entries for $\rho(t_{k,l},a_m)$ in the table within
Proposition~\ref{prop:ect} that are conjugates $y_i^w$ for some $i$ and
some $w \in (Y^\pm)^*$ can be replaced by $y_i$ for the purposes of
these calculations, and we shall henceforth denote such a conjugate of $y_i$
by $c(y_i)$.
So apart from the entry $y_n y_{n+1}^{-1}$ in the case $k=1, l<m=n$ (which will
be replaced by $z_n - z_{n+1}$), these are all words of length at most $1$.

We find that all the relators of $H$ within the set $S_2 = \{
	\rho(\epsilon,\varphi(y_i))y_i^{-1}: 1 \leq i \leq n+1$\}
are empty. For example, where $i=n+1$,
$\varphi(y_{n+1}) = a_2a_1\prod_{i=3}^n a_i$, and we have
$t_{1,2} a_2 = t_{1,3}$, $t_{1,3}a_1 = t_{2,3}$, $t_{2,l}a_l = t_{2,l+1}$
for $3 \leq l \leq n-1$, and $t_{2,n} a_n = y_{n+1} t_{1,2}$, so the
resulting relator of $H$ is $y_{n+1}y_{n+1}^{-1}$, freely reducing to the empty word.

So now we consider the relators of $H$ within the set
$S_1=\{\rho(t_{k,l},w) : t_{k,l} \in T, w \in R\}$.
We find that each relator of $H/[H,H]$ is derived several times, for various
values of $t_{k,l},w$.

The calculations of the abelianisations of relations of $H$ that are derived
from relations $a_ia_{i+1}a_i = a_{i+1}a_ia_{i+1}$ ($1 \leq i < n$),
$a_na_1a_n = a_1a_na_1$, or $a_ia_j = a_ja_i$ ($ 1 < i < j-1 \leq n-1$) of $G$
are routine. So we shall work out a couple of examples, and otherwise list
the results.

We derive relations $z_j=z_{j+1}$  for $H/[H,H]$ for all $i$ with
$2 \leq j \leq n-3$ using the braid relations
$w=a_ia_{i+1}a_i = a_{i+1}a_ia_{i+1}$ of $G$ with $1 \leq i \leq n-1$.

More precisely, 
exactly what is derived from the braid relation above depends on the 
associated transversal element $t=t_{k,l}$.
Suppose first that $i < n-1$. If $l<i$ then we derive $z_{i-1}=z_i$, 
while if $k<i<l-3$ then we derive $z_i=z_{i+1}$,
and if $i+1<k-1$ we derive $z_{i+1}=z_{i+2}$, but
all other cases give an empty relation.
Otherwise suppose that $i=n$. Then if $(k,l)=(1,n-1)$ we derive $z_{n-1}=z_n$,
while if $1<k<n-1$, $l=n$ we derive $z_{n-1}=z_n$, but 
all other cases give an empty relation.

The braid relation $a_na_1a_n = a_1a_na_1$ is more complicated, and
yields the empty relation for $(k,l) = (1,2)$ and $(1,n)$;
$z_n+z_2 = z_{n-2} + z_{n-1}$ when $1=k$, $3 \leq l < n$;
$z_{n-2}=z_2$ when $2 \leq k < l <n$;  $z_n=z_{n-1}$ when $(k,l) = (2,n)$;
and $z_{n-2}=z_2$ when $2 < k <l=n$.

For example, for the relation $a_na_1a_n = a_1a_na_1$ with $1=k$, $3\leq l < n$,
we get $t_{1,l}a_n = y_ny_{n+1}^{-1}t_{l,n}$, $t_{l,n}a_1 = y_2t_{l,n}$,
$t_{l,n}a_n = y_{n+1}t_{1,l}$; and
$t_{1,l}a_1 = t_{2,l}$, $t_{2,l}a_n = c(y_{n-2})t_{2,l}$,
$t_{2,l}a_1 = c(y_{n-1})y_{n-1}t_{1,l}$ which, as claimed, yields
$z_n+z_2 = z_{n-2} + z_{n-1}$ on abelianisation.

For the relations $a_ia_j = a_ja_i$ ($1 \leq i < j-1 \leq n-1$, $(i,j) \ne (1,n)$)
of $G$, there is a large but finite number of different cases to consider
depending on the relative values $i,j,k,l$. We claim first that they all yield
the empty relation of $H/[H,H]$ when $j < n$.

To see this, assume that $j<n$, and note first that, if
$k,l \not\in \{i,i+1,j,j+1\}$, then
$\overline{t_{k,l}a_i} = \overline{t_{k,l}a_j} = t_{k,l}$, and so the
resulting relation $\rho(t_{k,l},w_1) = \rho(t_{k,l},w_2)$ is just
$\rho(t_{k,l},a_i)\rho(t_{k,l},a_j) = \rho(t_{k,l},a_j)\rho(t_{k,l},a_i)$, which is trivial
on abelianisation.

Suppose next, that exactly one of $k,l$ is in $\{i,i+1,j,j+1\}$, say
$k \in \{i,i+1\}$ (the other three cases are similar).
Then $\overline{t_{k,l}a_j} = t_{k,l}$ and
$\overline{t_{k,l}a_i} = t_{k',l}$, with $k' = k \pm 1$, and then
$\overline{t_{k',l}a_j} = t_{k',l}$. We find from the table in Proposition~\ref{prop:ect}
that $\rho(t_{k,l},a_j) = \rho(t_{k',l},a_j) = y_{j-1}$ or $y_j$ (depending on
whether $l<j$ or $l>j$), so again the resulting relation of $H/[H,H]$
is trivial.

Finally (still assuming that $j<n$), suppose that $k,l \in \{i,i+1,j,j+1\}$,
so $k \in \{i,i+1\}$ and $l \in \{j,j+1\}$.
Then $\overline{t_{k,l}a_i} = t_{k',l}$, with $k' = k \pm 1$,
$\overline{t_{k,l}a_j} = t_{k,l'}$, with $l' = l \pm 1$, and
$\overline{t_{k,l'}a_i} = \overline{t_{k',l}a_j}=t_{k',l'}$.
Again, by using the table of Proposition~\ref{prop:ect}, we find that in each of the four
possible cases for $k$ and $l$, we get the trivial relation of $H/[H,H]$.

When $j=n$ we find, by similar calculations using the
table of Proposition~\ref{prop:ect}, that
we get the relation $z_i = z_{i+1}$ of $H/[H,H]$ when
$2 \leq i \leq n-2$ and ($k=1,i<l-1<n-1$ or $l=n, i<k-1$);
and the relation $z_{i-1} = z_i$ when $3 \leq i \leq n-1$ and
($k=1, i>l$ or $l=n, i>k$).

For example, if $2 \leq i \leq n-3$, $j=n$, $k=1$, and $i<l-1<n-1$, then we have
$t_{1,l}a_i = y_it_{1,l}$, $t_{1,l}a_n = y_ny_{n+1}^{-1}t_{l,n}$, and
$t_{l,n}a_i = y_{i+1}t_{l,n}$, so we get the relation $z_i=z_{i+1}$, as claimed.

Now all of these relations taken together reduce to
$z_2=z_3 = \cdots =z_{n-2}$ and $z_n=z_{n-1}$,
and hence $H/[H,H]$ is free abelian of rank $4$ with free basis the images
of $z_1$, $z_2$, $z_n$, $z_{n+1}$. This proves the first part of
Proposition~\ref{prop:aqhi}.

\subsection{The cycle and twisted cycle commutator relations of
$H_i/[H_i,H_i]$}
\label{sec:tabHrels}

The relations of $H_0$ and $H_i$ for $1 \leq t \leq n-2$ consist of all of the
relations of $H$ together with those that are derived from the cycle and
twisted cycle commutator relations of $G_0$ and $G_t$. These twisted
cycle commutator relations have the form $a_1 w = w a_1$, where
$w = v a_n v^{-1}$, and $v = a_2a_3 \cdots a_{n-1}$ for the cycle commutator
$\cc{a_1}{a_2}{\ldots}{a_n}$,
and $v  = a_2^{-1} a_3^{-1} \cdots a_{t+1}^{-1} a_{t+2} \cdots a_{n-1}$
for the twisted cycle commutator $\tc{a_1}{a_2}{\ldots}{a_n}_t$.

Since $\sigma(a_1)$ and $\sigma(w)$ both fix $\{1,2\}$ in all of the groups
$H_i$, we see that the relations $\rho(t_{k,l},a_1w) = \rho(t_{k,l},wa_1)$ of $H_i$
are trivial when $k>2$ and when $\{k,l\} = \{1,2\}$. 
So we need only consider the cases when $k=1$ or $2$ and $l>2$.

In these cases, we have $t_{1,l} a_1 = t_{2,l}$ and
$t_{2,l} a_1 = c(y_{n-1}) t_{1,l}$, whereas $t_{1,l} w = \alpha_{1l} t_{2,l}$ 
and $t_{2,l} w = \alpha_{2l} t_{1,l}$ for some words $\alpha_{1l}, \alpha_{2l}
\in (Y^{\pm})^*$ (which may depend also on which group $H_i$ we are
considering). Then the relations of $H_i$ in the cases $i=1$ and $2$ and $l>2$
are $\alpha_{2l} = \alpha_{1l} y_{n-1}$ and $y_{n-1} \alpha_{1l} = \alpha_{2l}$,
respectively, which have the same abelianisation, and hence we need only
consider the cases $k=1$, $3 \leq l \leq n$.

Before calculating the ensuing relations of $H_i/[H_i,H_i]$, it is helpful
to calculate the images of $\{1,l\}$ and of $\{2,l\}$ under the image
$\sigma(v)$ of the word $v$ under $\sigma$. Since
$\sigma(a_m) = \sigma(a_m^{-1})$ for all $i$, these are the same
in all of the groups $H_i$.

We have $\{1,l\}^{\sigma(a_m)} = \{1,l\}$ for $2 \leq m \leq l-2$,
$\{1,l\}^{\sigma(a_{l-1})} = \{1,l-1\}$, and
$\{1,l-1\}^{\sigma(a_m)} = \{1,l-1\}$ for $l \leq m \leq n-1$, so  
$\{1,l\}^{\sigma(v)} = \{1,l-1\}$.

For the image of $\{2,l\}$, we have
$\{m,l\}^{\sigma(a_m)} = \{m+1,l\}$ for $2 \leq m \leq l-2$,
$\{l-1,l\}^{\sigma(a_{l-1})} = \{l-1,l\}$, and
$\{l-1,m\}^{\sigma(a_m)} = \{l-1,m+1\}$ for $l \leq m \leq n-1$, so  
$\{2,l\}^{\sigma(v)} = \{l-1,n\}$.

In the case of $H_0$, the quotient by the normal closure of the cycle
commutator, we have
$t_{1,l} v = \beta_{1l} t_{1,l-1}$ and $t_{2,l} v = \beta_{2l} t_{l-1,n}$
for some words $\beta_{1l}, \beta_{2l} \in (Y^{\pm})^*$, and also
$t_{1,l-1} a_n = y_n y_{n+1}^{-1} t_{l-1,n}$ and 
$t_{l-1,n} a_n = y_{n+1} t_{1,k-1}$ so, denoting the images of
$\beta_{1l}, \beta_{2l}$ in $H_0/[H_0,H_0]$ by $\gamma_{1l}, \gamma_{2l}$,
the resulting relations of $H_0/[H_0,H_0]$ are
$\gamma_{2l} + z_{n+1} - \gamma_{1l} = \gamma_{1l} + z_n - z_{n+1} -
\gamma_{2l} + z_{n-1}$ or, equivalently,
\[2\gamma_{1l}-2\gamma_{2l} + z_{n-1} + z_n - 2z_{n+1} = 0.\]

A routine calculation shows that, for all $l$ with $3 \leq l \leq n$, we have
$\gamma_{1l} = \sum_{i=2}^{n-1} z_i$ and $\gamma_{2l} = z_1$, so from the
cycle commutator we get the single extra relation
\[-2z_1 + 2\sum_{i=2}^{n-2} z_i + 3z_{n-1} + z_n - 2z_{n+1}.\]

We saw above that $z_n$ and $z_{n-1}$ have equal images in $H/[H,H]$, so
\linebreak $H_0/[H_0,H_0] \cong \frac{\Z}{2\Z} \oplus \Z^3$. 

In the case of the quotients $H_t$ of $H$ by the normal closures of the
twisted cycle commutators $\tc{a_1}{a_2}{\ldots}{a_n}_t$, the corresponding
words in $(Y^{\pm})^*$ and their images in $H_t/[H_t,H_t]$ depend also on $t$,
and so we denote them by $\beta_{ilt}$ and $\gamma_{ilt}$ for $i=1,2$.
In this case, the corresponding calculation shows that we get two extra
relations when $1 \leq t \leq t-3$. (We observed at the beginning of this section
that $G_0 \cong G_{t-2}$ and hence $H_0 \cong H_{t-2}$.)

When $t+1 \leq l-2$ we have
\[ \gamma_{1lt} = -\sum_{i=2}^{t+1} z_i + \sum_{i=t+2}^{n-2} z_i + z_{n-1},\quad
\gamma_{2lt} = -tz_{n-1} + z_1, \] 
and when $t+1 \geq l-1$ we have
\[ \gamma_{1lt} = -\sum_{i=2}^{t} z_i + \sum_{i=t+1}^{n-2} z_i,\quad
\gamma_{2lt} = -(t-1)z_{n-1} - z_1. \] 

Here is some  more detail for the case $t+1 \geq l-1$.
For $2 \leq i \leq l-2$ we get $t_{1,l}a_i^{-1}= y_i^{-1}t_{1,l}$,
then $t_{1,l}a_{l-1}^{-1} = t_{1,l-1}$, then for $l \leq i \leq t+1$ we have
$t_{1,l-1}a_i^{-1}= y_{i-1}^{-1}t_{1,l-1}$, and finally
for $t+2 \leq i \leq n-1$ we have $t_{1,l-1}a_i= y_{i-1}t_{1,l-1}$, which results
in the claimed value of $\gamma_{1lt}$.

For $2 \leq i \leq l-2$ we get $t_{i,l}a_i^{-1}= c(y_{n-1})^{-1}t_{i+1,l}$, then
$t_{l-1,l} a_{l-1}^{-1} =y_1^{-1}t_{l-1,l}$, then for $l \leq i \leq t+1$ we have
$t_{l-1,i}t_i^{-1} = c(y_{n-1})^{-1} t_{l-1,i+1}$, and finally
for $t+2 \leq i \leq n-1$ we have $t_{l-1,i}a_i= t_{l-1,i+1}$,
so this results in the claimed value for $\gamma_{2lt}$.

The above equations result in the two relations
\begin{eqnarray*}
-2\sum_{i=1}^{t+1}z_i+2\sum_{i=t+2}^{n-2}z_i + (2t+3)z_{n-1}+z_n-2z_{n+1}&=&0,\\
2z_1-2\sum_{i=2}^{t}z_i+2\sum_{i=t+1}^{n-2}z_i+(2t-1)z_{n-1}+z_n-2z_{n+1}&=&0
\end{eqnarray*}
of $H_t/[H_t,H_t]$.
Subtracting the first of these from the second yields
$4z_1+4z_{t+1} -4z_{n-1}=0$ and, again using the fact that $z_{n-1}$ and $z_n$
have equal images in $H/[H,H]$, we see that 
 $H_t/[H_t,H_t] \cong \frac{\Z}{2\Z} \oplus \frac{\Z}{4\Z} \oplus \Z^2$,
for $1 \leq t \leq n$, which completes the proof of Proposition~\ref{prop:aqhi}.

\bibliographystyle{plain}
\bibliography{IsomNonisom.bib}

\end{document}